\documentclass[number,dvips]{arxbj}

\aid{0}
\volume{17}
\issue{1}
\pubyear{2011}
\firstpage{320}
\lastpage{346}
\doi{10.3150/10-BEJ270}

\makeatletter

\newtheorem{theorem}{Theorem}[section]

\newtheorem{lemma}{Lemma}[section]

\newtheorem{proposition}{Proposition}[section]

\newremark{remark}{Remark}[section]
\newproclaim{assumption}{Assumption}[section]
\newremark{example}{Example}[section]
\makeatother

\begin{document}
\begin{frontmatter}

\title{Mixing properties of ARCH and time-varying ARCH processes}
\runtitle{Mixing time-varying ARCH processes}

\begin{aug}
\author[a]{\fnms{Piotr} \snm{Fryzlewicz}\thanksref{a}\ead[label=e1]{p.fryzlewicz@lse.ac.uk}\corref{}}
\and
\author[b]{\fnms{Suhasini} \snm{Subba Rao}\thanksref{b}\ead[label=e2]{suhasini@stat.tamu.edu}}
\runauthor{P. Fryzlewicz and S. Subba Rao}
\address[a]{Department of Statistics,
London School of Economics,
Houghton Street,
London WC2A 2AE,
United Kingdom. \printead{e1}}
\address[b]{Department of Statistics,
Texas A\&M University,
3143 TAMU,
College Station, TX 77843-3143,
USA. \printead{e2}}
\end{aug}

\received{\smonth{11} \syear{2008}}
\revised{\smonth{2} \syear{2010}}

%
\begin{abstract}
There exist very few results on mixing for non-stationary processes.
However, mixing is often required in statistical inference for
non-stationary processes such as time-varying ARCH (tvARCH) models.
In this paper,
bounds for the mixing rates of a stochastic process are derived
in terms of the conditional densities of the process. These bounds
are used to obtain the $\alpha$, 2-mixing and $\beta$-mixing
rates of the non-stationary
time-varying $\operatorname{ARCH}(p)$ process and $\operatorname{ARCH}(\infty)$
process. It is shown that
the mixing rate of the time-varying $\operatorname{ARCH}(p)$ process
is geometric, whereas
the bound on the mixing rate of the $\operatorname{ARCH}(\infty)$ process
depends on the rate of decay of the $\operatorname{ARCH}(\infty)$ parameters.
We note that the methodology given in this paper is applicable to other
processes.
\end{abstract}

%

\end{frontmatter}

\section{Introduction}

Mixing is a measure of dependence between elements of a
random sequence that has a wide range of theoretical
applications (see \cite{pbra-07} and below).
One of the most popular mixing measures is
$\alpha$-mixing (also called \textit{strong mixing}),
where the $\alpha$-mixing rate of the non-stationary
stochastic process
$\{X_{t}\}$ is defined as a sequence of coefficients $\alpha(k)$ such that
\begin{eqnarray}\label{eq:strong-mix}
\alpha(k) =
\sup_{t\in\mathbb{Z}}\mathop{\sup_{H\in\sigma
(X_{t},X_{t-1},\ldots)}}_{G\in
\sigma(X_{t+k},X_{t+k+1},\ldots)}
|P(G\cap H) - P(G)P(H)|.
\end{eqnarray}
$\{X_{t}\}$ is called $\alpha$-mixing if $\alpha(k)\rightarrow0$ as
$k\rightarrow\infty$.
$\alpha$-mixing has several applications in statistical inference. For example,
if $\{\alpha(k)\}$ decays sufficiently fast to zero as $k\rightarrow
\infty$,
then, among other results, it is possible to show asymptotic normality of
sums of $\{X_{k}\}$ (see \cite{bdav-94}, Chapter 24),
as well as exponential inequalities for such sums
(see \cite{bbos-98}), asymptotic normality of kernel-based
nonparametric estimators
(see \cite{bbos-98}) and consistency of change point detection
schemes of nonlinear time series (see \cite{pfry-sub-08}).
The notion of 2-mixing is related to strong
mixing, but is a weaker condition as it measures the
dependence between two random variables and not the entire tails.
2-mixing is often used in statistical inference, for example, deriving
rates in nonparametric regression (see \cite{bbos-98}).
The 2-mixing rate can be used to
derive bounds for the covariance between functions of random variables, say
$\operatorname{cov}(g(X_{t}),g(X_{t+k}))$ (see \cite{pibr-62}),
which is usually not possible when only the correlation structure of
$\{X_{k}\}$ is known.
The 2-mixing rate of $\{X_{k}\}$ is defined as a sequence $\tilde
{\alpha}(k)$
which satisfies
\begin{eqnarray}\label{eq:2-mix}
\tilde{\alpha}(k)=
\sup_{t\in\mathbb{Z}}\mathop{\sup_{H\in\sigma(X_{t})}}_{G\in
\sigma(X_{t+k})}
|P(G\cap H) - P(G)P(H)|.
\end{eqnarray}
It is clear that $\tilde{\alpha}(k)\leq\alpha(k)$. A closely related
mixing measure, introduced in \cite{pvol-roz-59}
is $\beta$-mixing (also called \textit{absolutely
regular mixing}). The $\beta$-mixing rate of the stochastic process
$\{X_{t}\}$ is defined as a sequence of coefficients $\beta(k)$ such that
\begin{eqnarray}\label{eq:betamix}
\beta(k)=\sup_{t\in\mathbb{Z}}\mathop{\sup_{\{H_{j}\}\in\sigma(X_{t},X_{t-1},\ldots)}}_{\{G_{j}\}\in
\sigma(X_{t+k},X_{t+k+1},\ldots)}\sum_{i}\sum_{j}
|P(G_{i}\cap H_{j}) - P(G_{i})P(H_{j})|,
\end{eqnarray}
where $\{G_{i}\}$ and $\{H_{j}\}$ are finite partitions of
the sample space $\Omega$.
$\{X_{t}\}$ is called $\beta$-mixing if $\beta(k)\rightarrow0$ as
$k\rightarrow\infty$.
It can be seen that this measure is slightly stronger than
$\alpha$-mixing (since an upper bound for $\beta(k)$ immediately gives
a bound for $\alpha(k)$ due to the fact that $\beta(k)\geq\alpha(k)$).

Despite the versatility of mixing, its main drawback is that,
in general, it is difficult to derive bounds for
$\alpha(k)$, $\tilde{\alpha}(k)$ and
$\beta(k)$. However, the mixing bounds of some processes
are known. Chanda \cite{pcha-74}, Gorodetskii \cite{pgor-77},
Athreya and  Pantula \cite{path-pan-86} and Pham and  Tran \cite{ppha-85}
show strong mixing of the $\operatorname{MA}(\infty)$ process. Feigin and Tweedie \cite{pfei-twe-85} and
Pham \cite{ppha-86} have shown geometric ergodicity of bilinear processes
(we note that stationary geometrically ergodic Markov chains
are geometrically $\alpha$-mixing, $2$-mixing and $\beta$-mixing;
see, e.g., \cite{pfra-zak-06}).
More recently,
Tjostheim \cite{ptjo-90} and Mokkadem \cite{pmok-90} have shown geometric ergodicity for
a general class of Markovian processes. The results in \cite{pmok-90}
have been applied in \cite{pbou-98}
to show geometric ergodicity of stationary
$\operatorname{ARCH}(p)$ and $\operatorname{GARCH}(p,q)$ processes, where $p$ and $q$ are finite integers.
Related results on mixing for $\operatorname{GARCH}(p,q)$ processes can be found in
\cite{pcar-che-02,plie-05,plin-08,psor-06} (for
an excellent review) and \cite{pfra-zak-06,pmei-08} (where mixing of `nonlinear' $\operatorname{GARCH}(p,q)$ processes
is also considered).
Most of these these results are proved by verifying the Meyn--Tweedie
conditions (see \cite{pfei-twe-85} and \cite{bmey-twe-93}) and,
as mentioned above, are derived under the premise
that the process is stationary (or asymptotically stationary) and Markovian.
Clearly, if a process is non-stationary, then
the aforementioned results do not hold.
Therefore, for nonstationary processes, an alternative method to prove
mixing is required.

The main aim of this paper is to derive a bound for (\ref{eq:strong-mix}),
(\ref{eq:2-mix}) and (\ref{eq:betamix})
in terms of the densities of the process plus an additional
term, which is an extremal probability. These bounds can be applied to various
processes. In this paper, we will focus on ARCH-type processes and
use the bounds to derive mixing rates for time-varying $\operatorname{ARCH}(p)$
(tvARCH) and
$\operatorname{ARCH}(\infty)$ processes. The ARCH family of processes is widely
used in finance
to model the evolution of returns on financial instruments; we refer
the reader
to the review article of \cite{pgir-05} for a comprehensive overview of
mathematical
properties of ARCH processes and a list of further references.
It is worth mentioning that H{\"o}rmann \cite{hormann2008} and
Berkes \textit{et al.} \cite{berkeshormannschauer2009} have considered
a different type of dependence, namely a version of the
$m$-dependence moment measure, for ARCH-type processes.
The stationary $\operatorname{GARCH}(p,q)$ model tends to be the benchmark financial model.
However, in certain situations, it may not be the most appropriate
model. For example,
it cannot adequately explain the long memory seen in the data or change
according to shifts in the world economy. Therefore,
attention has recently been paid to tvARCH models (see, e.g.,
\cite{pdah-sub-06,pfry-sub-07,pfry-sub-08,pmik-sta-03}) and
$\operatorname{ARCH}(\infty)$ models (see \cite{pgir-00,pgir-rob-01,prob-91,psub-06}).
The derivations of the sampling properties of
some of the aforementioned papers rely on quite sophisticated assumptions
on the dependence structure, in particular, on their mixing properties.

We will show that, due to the $p$-Markovian nature of the time-varying
$\operatorname{ARCH}(p)$ process,
the $\alpha$-mixing, 2-mixing and $\beta$-mixing bounds have the same
geometric rate.
The story is different for $\operatorname{ARCH}(\infty)$ processes, where the mixing
rates can be
different and vary according to the rate of decay of the parameters.
An advantage of the
approach presented in this paper is that these methods can readily be
used to
establish mixing rates of
several time series models. This is especially useful in time series analysis,
for example, change point detection schemes for nonlinear time series, where
strong mixing of the underlying process is often required.
The price we pay for the flexibility of our approach is that the
assumptions under which we work are slightly stronger than the standard
assumptions required to prove geometric mixing of the stationary GARCH
process. However, the conditions
do not rely on proving irreducibility (which is usually required when
showing geometric
ergodicity) of the underlying process, which can be difficult to verify.

In Section \ref{sec:mix}, we derive a bound for the mixing rate of general
stochastic processes, in terms
of the differences of conditional densities. In Section \ref{sec:tvARCH},
we derive mixing bounds for time-varying $\operatorname{ARCH}(p)$ processes (where $p$
is finite).
In Section \ref{sec:ARCHinfty}, we derive mixing bounds for
$\operatorname{ARCH}(\infty)$
processes. Proofs which are not
in the main body of the paper can be found in the \hyperref[app]{Appendix} and the
accompanying technical report, available at
\url{http://stats.lse.ac.uk/fryzlewicz/mixing/tvARCH_mixing.pdf}.

\section{Some mixing inequalities for general processes}\label{sec:mix}

\subsection{Notation}

For $k>0$, let
$\underline{X}_{\,t}^{t-k} = (X_{t},\ldots,X_{t-k})$; if $k\leq0$,
then $\underline{X}_{\,t}^{t-k} = 0$.
Let $\underline{y}{}_{s} = (y_{s},\ldots,y_{0})$.
Let $\|\cdot\|$ denote the $\ell_{1}$-norm.
Let $\Omega$ denote the sample space.
The $\sigma$-algebra generated by $X_{t},\ldots,X_{t+r}$ is
denoted $\mathcal{F}_{t+r}^{t} =\sigma(X_{t},\ldots,X_{t+r})$.

\subsection{Some mixing inequalities}

Let us suppose that $\{X_{t}\}$ is an arbitrary stochastic process. In this
section, we derive some bounds for $\alpha(k)$, $\tilde{\alpha}(k)$ and
$\beta(k)$. To do this, we will consider bounds for
\begin{eqnarray*}
\sup_{H\in\mathcal{F}_{t}^{t-r_{1}},
G\in\mathcal{F}_{t+k+r_{2}}^{t+k}}
|P(G\cap H) - P(G)P(H)| \quad\mbox{and }
\\
\sup_{\{H_{j}\}\in\mathcal{F}_{t}^{t-r_{1}},
\{G_{i}\}\in
\mathcal{F}_{t+k+r_{2}}^{t+k}}\sum_{i,j}
|P(G_{i}\cap H_{j}) - P(G_{i})P(H_{j})|,
\end{eqnarray*}
where $r_{1}, r_{2} \geq0$ and $\{G_{i}\}$ and
$\{H_{i}\}$ are partitions of $\Omega$.
In the proposition below, we give a bound for the mixing
rate in terms of conditional densities. Similar bounds for
linear processes have been derived in \cite{pcha-74} and \cite{pgor-77}
(see also \cite{bdav-94}, Chapter 14). However, the bounds in
Proposition~\ref{propgeneral-bound} apply to any
stochastic process and it is this generality that allows us to use the result
in later
sections, where we derive mixing rates for ARCH-type processes.

\begin{proposition}\label{propgeneral-bound}
Let us suppose that the conditional density of
$\underline{X}_{\,t+k+r_{2}}^{t+k}$ given
$\underline{X}_{\,t}^{t-r_{1}}$ exists and denote it as
$f_{\underline{X}_{\,t+k+r_{2}}^{t+k}|\underline{X}_{\,t}^{t-r_{1}}}$.
For $\underline{\eta} =
(\eta_{0},\ldots,\eta_{r_{1}})\in(\mathbb{R}^{+})^{r_{1}+1}$,
define the set
\begin{eqnarray}\label{eq:defE}
E = \{\omega;  \underline{X}_{\,t}^{t-r_{1}}(\omega)\in\mathcal{E}
\},\qquad
\mbox{where }
\mathcal{E} = \{(\nu_{0},\ldots,\nu_{r_{1}});
\mbox{ for all } |\nu_{j}| \leq\eta_{j}\}.
\end{eqnarray}
For all $r_{1},r_{2}\geq0$ and $\underline{\eta}$, we then have
\begin{eqnarray}\label{eq:mix-bd}
&&\sup_{H\in\mathcal{F}_{t}^{t-r_{1}}, G\in
\mathcal{F}_{t+k+r_{2}}^{t+k}}
|P(G\cap H) - P(G)P(H)| \nonumber
\\[-8pt]\\[-8pt]
&&\quad \leq
2 \sup_{\underline{x}\in\mathcal{E}}\int_{\mathbb{R}^{r_{2}+1}}
\bigl|f_{\underline{X}_{\,t+k+r_{2}}^{t+k}|\underline{X}_{\,t}^{t-r_{1}}}(
\underline{y}|\underline{x})-
f_{\underline{X}_{\,t+k+r_{2}}^{t+k}|\underline{X}_{\,t}^{t-r_{1}}}
(\underline{y}|0)\bigr|\,
\mathrm{d}\underline{y} +
4P(E^{c})\nonumber
\end{eqnarray}
and
\begin{eqnarray}\label{eq:mix-bdii}
&&\sup_{\{H_{j}\}\in\mathcal{F}_{t}^{t-r_{1}},\{G_{j}\}\in
\mathcal{F}_{t+k+r_{2}}^{t+k}}\sum_{i,j}
|P(G_{i}\cap H_{j}) - P(G_{i})P(H_{j})| \nonumber
\\[-8pt]\\[-8pt]
&&\quad \leq
2 \int_{\mathbb{R}^{r_{2}+1}}
\sup_{\underline{x}\in\mathcal{E}}
\bigl|f_{\underline{X}_{\,t+k+r_{2}}^{t+k}|\underline{X}_{\,t}^{t-r_{1}}}(
\underline{y}|\underline{x})-
f_{\underline{X}_{\,t+k+r_{2}}^{t+k}|\underline{X}_{\,t}^{t-r_{1}}}
(\underline{y}|0)\bigr|
\,\mathrm{d}\underline{y} + 4P(E^{c}),\nonumber
\end{eqnarray}
where $\{G_{i}\}$ and $\{H_{j}\}$ are finite partitions of $\Omega$.
Letting $\underline{W}^{t+1}_{t+k-1}$ be a random vector that is
independent of $\underline{X}_{\,t}^{t-r_{1}}$,
we then have
\begin{eqnarray}\label{eq:mix-bd2}
&&\sup_{H\in\mathcal{F}_{t}^{t-r_{1}}, G\in
\mathcal{F}_{t+k+r_{2}}^{t+k}}
|P(G\cap H) - P(G)P(H)| \nonumber
\\[-8pt]\\[-8pt]
&&\quad \leq 2 \sum_{s=0}^{r_{2}} \sup_{\underline{x}\in\mathcal{E}}
\mathbb{E}_{\underline{W}}\biggl(\sup_{\underline{y}{}_{s-1}\in
\mathbb{R}^{s}}
\int_{\mathbb{R}}
\mathcal{D}_{s,k,t}(y_{s}|\underline{y}{}_{s-1},\underline{W},
\underline{x})
\,\mathrm{d} y_{s}\biggr) + 4P(E^{c})\nonumber
\end{eqnarray}
 and
\begin{eqnarray}\label{eq:mix-bd2ii}
&& \sup_{\{H_{j}\}\in\mathcal{F}_{t}^{t-r_{1}},\{G_{j}\}\in
\mathcal{F}_{t+k+r_{2}}^{t+k}}\sum_{i,j}
|P(G_{i}\cap H_{j}) - P(G_{i})P(H_{j})| \nonumber
\\[-8pt]\\[-8pt]
&&\quad \leq
2 \sum_{s=0}^{r_{2}}\mathbb{E}_{\underline{W}}
\biggl(\sup_{\underline{y}{}_{s-1}\in\mathbb{R}^{s}}
\int_{\mathbb{R}}\sup_{\underline{x}\in\mathcal{E}}
\mathcal{D}_{s,k,t}(y_{s}|\underline{y}{}_{s-1},\underline{W},
\underline{x})\,\mathrm{d} y_{s}\biggr) + 4P(E^{c}),\nonumber
\end{eqnarray}
where $\mathbb{E}_{\underline{W}}(g(\underline{W})) = \int
g(\underline{w})
f_{\underline{W}}(\underline{w})\,\mathrm{d}\underline{w}$,
$f_{\underline{W}}(\underline{w})$ is the density of $\underline{w}$,
$\mathcal{D}_{0,k,t}(y_{0}|
\underline{y}{}_{-1},\underline{w},\underline{x}) =
|f_{s,k,t}(y_{s}|\underline{w},\underline{x})-
f_{s,k,t}(y_{s}|\underline{w},0)|$ and, for $s\geq1$,
\begin{eqnarray}
\label{eq:G}
\mathcal{D}_{s,k,t}(y_{s}|
\underline{y}{}_{s-1},\underline{w},\underline{x})
=\bigl|f_{s,k,t}
(y_{s}|\underline{y}{}_{s-1},\underline{w},\underline{x}) -
f_{s,k,t}(y_{s}|\underline{y}{}_{s-1},\underline{w},0)\bigr|
\end{eqnarray}
with the conditional density of
$X_{t+k}$ given $(\underline{W}_{t+k-1}^{t+1},\underline{X}_{\,t}^{t-r_{1}})$
denoted $f_{0,k,t}$,
the conditional density of
$X_{t+k+s}$ given $(\underline{X}^{\,t+k}_{t+k+s-1},
\underline{W}_{t+k-1}^{t+1},\underline{X}_{\,t}^{t-r_{1}})$
denoted $f_{s,k,t}$,
$\underline{x} = (x_{0},\ldots,x_{-r_{2}})$
and $\underline{w} = (w_{k},\ldots,w_{1})$.
\end{proposition}

\begin{pf}
This can be found in Appendix \hyperref[sec:app1]{A.1}.
\end{pf}

Since the above bounds hold for all vectors $\underline{\eta}\in
(\mathbb{R}^{+})^{r_{1}+1}$
(note that $\underline{\eta}$ defines the set $E$; see (\ref{eq:defE})),
by choosing the $\underline{\eta}$ which balances the integral and
$P(E^{c})$, we obtain an upper bound for the mixing rate.

The main application of the inequality in (\ref{eq:mix-bd2})
is to processes which
are `driven' by the innovations (e.g., linear and ARCH-type processes).
If $\underline{W}_{t+k-1}^{t+1}$ is the
innovation process, it can often be shown that
the conditional density of $X_{t+k+s}$ given
$(\underline{X}^{\,t+k}_{t+k+s-1},
\underline{W}_{t+k-1}^{t+1},\underline{X}_{\,t}^{t-r_{1}})$ can be
written as
a function of the innovation density. Deriving the density
of $X_{t+k+s}$ given $(\underline{X}^{t+k}_{\,t+k+s-1},
\underline{W}_{t+k-1}^{t+1},\underline{X}_{\,t}^{t-r_{1}})$
is not a trivial task, but it is often possible.
In the subsequent sections, we will apply Proposition~\ref{propgeneral-bound}
to obtain bounds for the mixing rates.

The proof of Proposition \ref{propgeneral-bound}
can be found in the \hyperref[app]{Appendix}, but we give a brief outline of it here. Let
\begin{eqnarray}\label{eq:HG}
 H = \{\omega; \underline{X}_{\,t}^{t-r_{1}}(\omega) \in\mathcal
{H}\}, \qquad
G = \{\omega; \underline{X}_{\,t+k+r_{2}}^{t+k}(\omega) \in\mathcal
{G}\}.
\end{eqnarray}
It is straightforward to show that
$|P(G\cap H) - P(G)P(H)| \leq
|P(G\cap H\cap E) - P(G\cap E)P(H)| + 2P(E^{c})$.
The advantage of this decomposition is that when we restrict
$\underline{X}_{\,t}^{t-r_{1}}$ to the set $\mathcal{E}$ (i.e., not
large values of $\underline{X}_{\,t}^{t-r_{1}}$),
we can obtain a bound for
$|P(G\cap H\cap E) - P(G\cap E)P(H)|$. More precisely, by using
the inequality
\begin{eqnarray*}
\inf_{\underline{x}\in\mathcal{E}}
P(G|\underline{X}_{\,t}^{t-r_{1}}=\underline{x})
P(H\cap E) \leq
P(G\cap H\cap E)
\leq \sup_{\underline{x}\in\mathcal{E}}
P(G|\underline{X}_{\,t}^{t-r_{1}}=\underline{x})
P(H\cap E),
\end{eqnarray*}
%
we can derive upper and lower bounds for
$P(G\cap H\cap E) - P(G\cap E)P(H)$ which depend
only on $E$ and not $H$ and $G$, and thus obtain the
bounds in Proposition \ref{propgeneral-bound}.

It is worth mentioning that by using (\ref{eq:mix-bd2}), one can establish
mixing rates for time-varying linear processes (such as the
tvMA$(\infty)$ process
considered in \cite{pdah-pol-06}). Using (\ref{eq:mix-bd2}) and
techniques similar to those used in Section \ref{sec:ARCHinfty},
mixing bounds can be
obtained for the tvMA$(\infty)$ process.

In the following
sections, we will derive the mixing rates for ARCH-type processes, where
one of the challenging aspects of the proof is establishing a bound
for the integral difference in~(\ref{eq:G}).

\section{Mixing for the time-varying $\operatorname{ARCH}(p)$ process}\label{sec:tvARCH}

\subsection{The tvARCH process}

In \cite{pfry-sub-07}, it is shown that the tvARCH process can be used
to explain
the commonly observed stylized facts in financial time series
(such as the empirical long memory).
A sequence of random variables $\{X_{t}\}$ is said to come from the
squares of a time-varying $\operatorname{ARCH}(p)$ process if it satisfies the representation
\begin{eqnarray}
\label{eq:tvARCH}
X_{t} = Z_{t}\Biggl(a_{0}(t) + \sum_{j=1}^{p}a_{j}(t)X_{t-j} \Biggr),
\end{eqnarray}
where $\{Z_{t}\}$ are independent, identically distributed
(i.i.d.) positive random variables, where $\mathbb{E}(Z_{t})=1$ and
$a_{j}(\cdot)$ are positive parameters.
It is worth comparing (\ref{eq:tvARCH}) with
the squared tvARCH process used in the statistical literature. Unlike the
squared tvARCH process considered in, for example, \cite{pdah-sub-06} and
\cite{pfry-sub-07}, we have not placed {\it any} smoothness
conditions on the time-varying parameters $\{a_{j}(\cdot)\}$.
The smoothness conditions assumed in \cite{pdah-sub-06} and
\cite{pfry-sub-07} are used in order to carry out parameter
estimation. However,
in this paper, we are dealing with mixing of the process, which does not
require such strong assumptions.
The assumptions that we require are stated below.
From now on, with a slight abuse of terminology, we will
call the squared tvARCH process simply the tvARCH process.

\begin{assumption}\label{assum}
\textup{(i)} For some $\delta> 0$,
$\sup_{t\in\mathbb{Z}}\sum_{j=1}^{p}a_{j}(t)\leq1-\delta$.\vspace*{-5pt}
\begin{longlist}[(iii)]
\item[(ii)] $\inf_{t\in\mathbb{Z}}
a_{0}(t)> 0 $ and
$\sup_{t\in\mathbb{Z}}a_{0}(t)<\infty$.
\item[(iii)] Let $f_{Z}$ denote the density of $Z_{t}$. For all
$a >0$, we have
$\int| f_{Z}(u) -
f_{Z}(u[1+a])|\,\mathrm{d}u \leq K a$
for some finite $K$ independent of $a$.
\item[(iv)] Let $f_{Z}$ denote the density of $Z_{t}$. For all
$a >0$, we have $\int\sup_{0 \leq\tau\leq a}| f_{Z}(u) -
f_{Z}(u[1+\tau])|\,\mathrm{d}u \leq K a$
for some finite $K$ independent of $a$.
\end{longlist}
\end{assumption}

We note that Assumption \ref{assum}(i)--(ii) guarantees that the ARCH
process has a Volterra expansion as a solution (see \cite{pdah-sub-06},
Section 5).
Assumption~\ref{assum}(iii)--(iv) is a type of Lipschitz condition on the
density function and is satisfied by various well-known
distributions, including the chi-squared distributions.
We now consider
a class of densities which satisfy Assumption~\ref{assum}(iii)--(iv).
Suppose that $f^{\prime}_{Z}$ is bounded, that
after some finite point $m$ the derivative $f^{\prime}$ declines
monotonically to
zero and satisfies $\int|y f_{Z}^{\prime}(y)|\,\mathrm{d}y < \infty$. In this case,
\begin{eqnarray*}
&& \int_{0}^{\infty}\sup_{0 \leq\tau\leq a}| f_{Z}(u) -
f_{Z}(u[1+\tau])|\,\mathrm{d}u
\\
&&\quad \leq \int_{0}^{m} \sup_{0 \leq\tau\leq a}
| f_{Z}(u) -
f_{Z}(u[1+\tau])|\,\mathrm{d}u +
\int^{\infty}_{m}
\sup_{0 \leq\tau\leq a}| f_{Z}(u) - f_{Z}(u[1+\tau])
|\,\mathrm{d}u
\\
&&\quad \leq a\biggl( m^{2}\sup_{u\in\mathbb{R}}
|f_{Z}^{\prime}(u)|+ \int^{\infty}_{m}u|f_{Z}^{\prime}(u)|\,\mathrm{d}u\biggr)
\leq Ka
\end{eqnarray*}
for some finite $K$ independent of $a$, hence Assumption \ref{assum}(iii)--(iv)
is satisfied.

We use Assumption \ref{assum}(i)--(iii) to obtain the strong
mixing rate (2-mixing and $\alpha$-mixing) of the tvARCH$(p)$ process, and
the slightly stronger conditions Assumption \ref{assum}(i)--(ii) and (iv)
to obtain the $\beta$-mixing rate of the tvARCH$(p)$
process. We mention that in the case that $\{X_{t}\}$ is a stationary,
ergodic time series, \cite{pfra-zak-06} have shown geometric ergodicity,
which they show implies $\beta$-mixing, under the weaker condition
that the
distribution function of $\{Z_{t}\}$ can have some discontinuities.


\subsection{The tvARCH$(p)$ process and the Volterra series expansion}

In this section, we derive a Volterra series expansion of the tvARCH
process (see also \cite{pgir-00}).
These results allow us to apply Proposition \ref{propgeneral-bound}
to the
tvARCH process.
We first note that the innovations $\underline{Z}_{t+k-1}^{t+1}$
and $\underline{X}^{t-p+1}_{\,t}$ are independent random vectors. Hence,
comparing with Proposition \ref{propgeneral-bound}, we are interested in
obtaining the conditional density
of $X_{t+k}$ given $\underline{Z}_{t+k-1}^{t+1}$ and
$\underline{X}^{t-p+1}_{t}$
(denoted $f_{0,k,t}$) and the conditional density
of $X_{t+k+s}$ given $\underline{X}_{\,t+k+s-1}^{t+k},
\underline{Z}_{t+k-1}^{t+1}$ and $\underline{X}^{t-p+1}_{\,t}$ (denoted
$f_{s,k,t}$).
We use these expressions to obtain a bound for $\mathcal{D}_{s,k,t}$
(defined in (\ref{eq:G})),
which we use to derive a bound for the mixing rate.
We now represent $\{X_{t}\}$ in terms of $\{Z_{t}\}$. To do this, we define\vspace*{-2pt}
\begin{eqnarray*}
A_{t} &=&
\pmatrix{
a_{1}(t)Z_{t} & a_{2}(t)Z_{t} & \ldots& a_{p}(t)Z_{t} \cr
1 & 0 & \ldots& 0 \cr
0 & 1 & \ldots& 0 \cr
\ldots & \ldots & \ddots& \vdots\cr
0 & 0 & 1 & 0 \cr
}
,\qquad\tilde{A}_{t}
=
\pmatrix{
a_{1}(t) & a_{2}(t) & \ldots& a_{p}(t) \cr
1 & 0 & \ldots& 0 \cr
0 & 1 & \ldots& 0 \cr
\ldots & \ldots & \ddots& \vdots\cr
0 & 0 & 1 & 0 \cr
}
,
\\[-2pt]
\underline{b}_{t}
&=& (a_{0}(t)Z_{t},0,\ldots,0)^{\prime} \quad\mbox{and} \quad
\underline{X}_{\,t}^{t-p+1} = (X_{t},X_{t-1},\ldots,X_{t-p+1})^{\prime}.\vspace*{-2pt}
\end{eqnarray*}
%
Using this notation, we have the relation
$\underline{X}_{\,t+k}^{t+k-p+1}=A_{t+k}\underline{X}_{\,t+k-1}^{t+k-p}
+ \underline{b}_{t+k}$.
We note that the vector representation of ARCH and GARCH processes has been
used in \cite{pbas-03,pbou-pic-92a,pstr-mik-06} in order to
obtain some probabilistic properties for ARCH-type processes.
Now iterating, the relation $k$ times (to get $\underline{X}_{\,t+k}^{t+k-p+1}$
in terms of $\underline{X}_{\,t}^{t-p+1}$), we have
\begin{eqnarray*}
\label{eq:underlineX}
\underline{X}_{\,t+k}^{t+k-p+1} = \underline{b}_{t+k} +
\sum_{r=0}^{k-2}\Biggl[\prod_{i=0}^{r-1}A_{t+k-i}\Biggr]
\underline{b}_{t+k-r-1} +
\Biggl[\prod_{i=0}^{k-1}A_{t+k-i}\Biggr]\underline{X}_{\,t}^{t-p+1},
\end{eqnarray*}
where we set $[\prod_{i=0}^{-1}A_{t+k-i}] = I_{p}$
($I_{p}$ denotes the $p\times p$-dimensional identity matrix).
We use this expansion below.

\begin{lemma}
\label{lemma:3.1}
Let us suppose that Assumption \textup{\ref{assum}(i}) is satisfied. For $s\geq
0$, we then have
\begin{eqnarray}\label{eq:tvARCHEXPa}
X_{t+k+s} &=& Z_{t+k+s}\{
\mathcal{P}_{s,k,t}(\underline{Z})
+ \mathcal{Q}_{s,k,t}(
\underline{X})\} ,
\end{eqnarray}
where $\underline{Z} = \underline{Z}_{t+k}^{t+1}$;
for $s=0$ and $n>t$, we have
$\mathcal{P}_{0,k,t}(\underline{Z})=
a_{0}(t+k) + [\tilde{A}_{t+k} \times\break \sum_{r=0}^{n-t-2}
\prod_{i=1}^{r}A_{t+k-i}b_{t+k-r-1}]_{1}$,
$\mathcal{Q}_{0,k,t}(\underline{X})
= [\tilde{A}_{t+k}\prod_{i=1}^{k-1}A_{t+k-i}
\underline{X}_{\,t}^{t-p+1}]_{1}$
($[\cdot]_{1}$ denotes the first element of a vector).

For $1\leq s\leq p$,
\begin{eqnarray}
\label{eq:Q}
\mathcal{P}_{s,k,t}(\underline{Z})
&=& a_{0}(t+k+s) + \sum_{i=1}^{s-1}a_{i}(t+k+s)X_{t+k+s-i}\nonumber
\\
&&{}+ \sum_{i=s}^{p}a_{i}(t+k+s)Z_{k+s-i}\nonumber
\\
&&{}\hspace*{23pt}\times\Biggl\{a_{0}(t+k+s-i)
\\
&&{}\hspace*{40pt}+
\Biggl[\tilde{A}_{t+k+s-i} \sum_{r=1}^{k+s-i}
\Biggl\{\prod_{d=0}^{r}A_{t+k+s-i-d}\Biggr\}b_{t+k+s-i-r}\Biggr]_{1}\Biggr\},\nonumber
\\
 \mathcal{Q}_{s,k,t}(\underline{Z},\underline{X}) &=&
\Biggl[\sum_{i=s}^{p}a_{i}(t+k+s)Z_{k+s-i}
\tilde{A}_{t+k+s-i}\Biggl\{\prod_{d=0}^{k+s-i}A_{t+k+s-i-d}
\underline{X}_{\,t}^{t-p+1}\Biggr\}\Biggr]_{1} \nonumber
\end{eqnarray}
and for $s> p$, we have
$\mathcal{P}_{s,k,t}(\underline{Z}) =
a_{0}(t+k+s) + \sum_{i=1}^{p}a_{i}(t+k+s)X_{t+k+s-i}$
and $\mathcal{Q}_{s,k,t}(\underline{Z},\underline{X})\equiv0$. We
note that
$\mathcal{P}_{s,k,t}(\cdot)$ and
$\mathcal{Q}_{s,k,t}(\cdot)$ are positive random variables
and for $s\geq1$, $\mathcal{P}_{s,k,t}(\cdot)$ is a function of
$\underline{X}_{\,t+k+s-1}^{t+k}$ (but this has been suppressed in the
notation).
\end{lemma}

\begin{pf}
This is found in Appendix \hyperref[sec:app2]{A.2}.
\end{pf}

By using (\ref{eq:tvARCHEXPa}), we now show that the conditional
density of
$X_{t+k+s}$ given $\underline{X}^{t+k}_{\,t+k+s-1},\break\underline
{Z}_{t+k-1}^{t+1}$ and
$\underline{X}_{\,t}^{t-p+1}$ is a function of the density of $Z_{t+k+s}$.
It is clear from (\ref{eq:tvARCHEXPa}) that
$Z_{t+k+s}$ can be expressed as
$Z_{t+k+s} = \frac{X_{t+k+s}}{\mathcal{P}_{s,k,t}(
\underline{Z}) + \mathcal{Q}_{s,k,t}(\underline{Z},
\underline{X})}$.
Therefore, it is straightforward to show that
\begin{eqnarray}\label{eq:fxz}
f_{s,k,t}(y_{s}|\underline{y}{}_{s-1},\underline{z},\underline{x})
=\frac{1}{\mathcal{P}_{s,k,t}(\underline{z})+
\mathcal{Q}_{s,k,t}(\underline{z},\underline{x})}
f_{Z}\biggl(\frac{y_{s}}{\mathcal{P}_{s,k,t}(
\underline{z})
+\mathcal{Q}_{s,k,t}(\underline{z},\underline{x})}\biggr).
\end{eqnarray}


\subsection{Strong mixing of the tvARCH$(p)$ process}\label{sec:density:archp}

The aim of this section is to prove geometric mixing of the tvARCH$(p)$
process without appealing to geometric ergodicity. Naturally, the
results in this section also apply to stationary $\operatorname{ARCH}(p)$ processes.

In the following lemma, we use Proposition \ref{propgeneral-bound}
to obtain bounds for the mixing rates.
It is worth mentioning that the techniques
used in the proof
below can be applied to other Markov processes.

\begin{lemma}\label{lemma:smaller-algebra}
Suppose that $\{X_{t}\}$ is a tvARCH process which satisfies (\ref{eq:tvARCH}).
For any
$\underline{\eta} = (\eta_{0},\ldots,\eta_{-p+1})\in
(\mathbb{R}^{+})^{p}$, we then have
\begin{eqnarray}\label{eq:diff-f}
&&\sup_{G\in\mathcal{F}^{t+k}_{\infty},H\in
\mathcal{F}^{-\infty}_{t}}
|P(G\cap H) - P(G)P(H)| \nonumber
\\
&&\quad \leq
2\sum_{s=0}^{p-1} \sup_{\underline{x}\in\mathcal{E}} \int
\mathbb{E}_{\underline{Z}}\biggl(\sup_{\underline{y}{}_{s-1}\in
\mathbb{R}^{s}}
\int_{\mathbb{R}}
\mathcal{D}_{s,k,t}(y_{s}|
\underline{y}{}_{s-1},\underline{Z},\underline{x})\,\mathrm{d}y_{s}\biggr)
\\
&&{}\qquad + 4\sum_{j=0}^{p-1}P(|X_{t-j}|\geq\eta_{-j+1}) \nonumber
\end{eqnarray}
and
\begin{eqnarray}\label{eq:diff-f2}
&& \sup_{\{H_{j}\}\in\mathcal{F}_{t}^{-\infty},\{G_{j}\}\in
\underline{F}_{\,\infty}^{t+k}}\sum_{i,j}
|P(G_{i}\cap H_{j}) - P(G_{i})P(H_{j})|\nonumber
\\
&&\quad \leq
2\sum_{s=0}^{p-1} \sup_{\underline{x}\in\mathcal{E}}
\mathbb{E}_{\underline{Z}}\biggl(\sup_{\underline{y}{}_{s-1}\in
\mathbb{R}^{s}}
\int_{\mathbb{R}}
\sup_{\underline{x}\in\mathcal{E}}
\mathcal{D}_{s,k,t}(y_{s}|
\underline{y}{}_{s-1},\underline{Z},\underline{x})\,\mathrm{d}y_{s}
\biggr)
\\
&&{}\qquad + 4\sum_{j=0}^{p-1}P(|X_{t-j}|\geq\eta_{-j+1}), \nonumber
\end{eqnarray}
where
$\underline{z} = (z_{1},\ldots,z_{k-1})$ and
$\{G_{i}\}$ and $\{H_{j}\}$ are partitions of $\Omega$ and
$\mathbb{E}_{\underline{Z}}(g(\underline{Z})) = \int g(\underline{z})
\times\break \prod_{i=1}^{k-1}f_{Z}(z_{i})\,\mathrm{d}z_{i}$.
\end{lemma}

\begin{pf} This can be found in Appendix \hyperref[sec:app2]{A.2}.
\end{pf}

To obtain a mixing rate for the tvARCH$(p)$ process, we need to bound
the integral in (\ref{eq:diff-f}), then
obtain the set $E$ which minimizes (\ref{eq:diff-f}).
We will start by bounding
$\mathcal{D}_{s,k,t}$,
which, we recall, is based on the conditional density
$f_{s,k,t}$ (defined in (\ref{eq:fxz})).

\begin{lemma}
\label{lemma:H}
Let $ \mathcal{D}_{s,k,t}$ and $\mathcal{Q}_{s,k,t}$
be defined as in (\ref{eq:G}) and (\ref{eq:Q}), respectively.
\begin{enumerate}[(ii)]
\item[(i)] Supposing that Assumption \textup{\ref{assum}(i)--(iii)} holds,
then for all
$\underline{x}\in(\mathbb{R}^{+})^{p}$, we have
\begin{eqnarray}
\label{eq:density-difference}
\sum_{s=0}^{p-1} \int\mathbb{E}_{\underline{Z}}\biggl(
\sup_{\underline{y}{}_{s-1}\in\mathbb{R}^{s}}
\int
\mathcal{D}_{s,k,t}(y_{s}|
\underline{y}{}_{s-1},\underline{Z},\underline{x})\,\mathrm{d}y_{s}
\biggr)&\leq& K\frac{\mathbb{E}[\mathcal{Q}_{s,k,t}(\underline{Z},
\underline{x})]}{\inf_{t\in\mathbb{Z}}a_{0}(t)}\nonumber
\\[-8pt]\\[-8pt]
&\leq&
K(1-\tilde{\delta})^{k}\|\underline{x}\|,\nonumber
\end{eqnarray}
where $K$ is a finite constant and $0<\tilde{\delta}\leq\delta<1$
($\delta$ is defined in Assumption~\textup{\ref{assum}(i)}).
\item[(ii)] Supposing that Assumption~\textup{\ref{assum}(i)--(ii)} and \textup{(vi)}
hold, then for any
set $\mathcal{E}$ (defined as in (\ref{eq:defE})), we have
\begin{eqnarray}
\label{eq:density-difference2}
\sum_{s=0}^{p-1} \mathbb{E}_{\underline{Z}}\biggl(
\sup_{\underline{y}{}_{s-1}\in
\mathbb{R}^{s}} \int\sup_{\underline{x}\in\mathcal{E}}
\mathcal{D}_{s,k,t}(y_{s}|
\underline{y}{}_{s-1},\underline{Z},\underline{x})\,\mathrm{d}y_{s}
\biggr) \leq
\sup_{\underline{x}\in\mathcal{E}}
K(1-\tilde{\delta})^{k}\|\underline{x}\|.
\end{eqnarray}
\end{enumerate}
\end{lemma}

\begin{pf} This can be found in Appendix \hyperref[sec:app2]{A.2}.
\end{pf}

We now use the lemmas above to show geometric mixing of the tvARCH process.

\begin{theorem}\label{thm:archp}
\textup{(i)} Supposing that Assumption~\textup{\ref{assum}(i)--(iii)} holds, then
\begin{eqnarray*}
\mathop{\sup_{G\in\sigma(\underline{X}^{t+k}_{\,\infty})}}_{H\in
\sigma(\underline{X}^{-\infty}_{\,t})}|P(G\cap H) -
P(G)P(H)| \leq K\alpha^{k}.
\end{eqnarray*}\vspace*{-8pt}
\begin{longlist}[(ii)]
\item[(ii)] Supposing that Assumption \textup{\ref{assum}(i)--(ii)} and \textup{(iv)}
hold, then
\begin{eqnarray*}
\mathop{\sup_{\{H_{j}\}\in\sigma(\underline{X}_{\,t}^{-\infty})}}_{\{
G_{j}\}\in
\sigma(\underline{X}_{\,\infty}^{t+k})}\sum_{i}\sum_{j}
|P(G_{i}\cap H_{j}) - P(G_{i})P(H_{j})|
\leq K\alpha^{k}
\end{eqnarray*}
\end{longlist}
for any $\sqrt{1-\delta} <\alpha<1$, where $K$ is a finite constant
independent of
$t$ and $k$.
\end{theorem}

\begin{pf} We will use (\ref{eq:diff-f}) to prove (i).
Equation (\ref{eq:density-difference}) gives a bound
for the integral difference in (\ref{eq:diff-f}); therefore, all that remains
is to bound the probabilities in (\ref{eq:diff-f}).
To do this, we first use Markov's inequality, to give
$\sum_{j=0}^{p-1}P(|X_{t-j}|\geq\eta_{-j}) \leq
\sum_{j=0}^{p-1}\mathbb{E}|X_{t-j}|\eta_{-j}^{-1}$. By using the Volterra
expansion of $X_{t}$ (see \cite{pdah-sub-06}, Section 5), it can be
shown that $\sup_{t\in\mathbb{Z}}\mathbb{E}|X_{t}|\leq
(\sup_{t\in\mathbb{Z}}a_{0}(t))/(1-\sup_{t\in\mathbb{Z}}\sum_{j=1}^{p}
a_{j}(t))$. Using these bounds and
substituting (\ref{eq:density-difference}) into
(\ref{eq:diff-f}) gives, for every $\underline{\eta}\in(\mathbb
{R}^{+})^{p}$,
the bound
\begin{eqnarray*}
\mathop{\sup_{G\in\sigma(\underline{X}^{t+k}_{\,\infty})}}_{H\in
\sigma(\underline{X}^{-\infty}_{\,t})}
|P(G\cap H) - P(G)P(H)|
\leq 2\frac{K(1-\tilde{\delta})^{k}\sum_{j=0}^{p-1}\eta_{-j}}{
\inf_{t\in\mathbb{Z}}a_{0}(t)}
+ 4K\sum_{j=0}^{p-1}\frac{1}{\eta_{-j}}.
\end{eqnarray*}
We observe that the right-hand side of the above is minimized when
$\eta_{-j} = (1-\tilde{\delta})^{k/2}$ (for
$0\leq j \leq(p-1)$), which gives the bound
\begin{eqnarray*}
\mathop{\sup_{H\in\sigma(\underline{X}^{-\infty}_{\,t})}}_{G\in
\sigma(\underline{X}_{\,\infty}^{t+k})}
|P(G\cap H) - P(G)P(H)| \leq K \sqrt{(1-\tilde{\delta})^{k}}.
\end{eqnarray*}
Since the above is true for any $0 < \tilde{\delta} < \delta$,
(ii) is true for any $\alpha$ which satisfies
$\sqrt{1-\delta} < \alpha< 1$, thus giving the result.

To prove (ii), we use an identical argument, but using the bound in
(\ref{eq:density-difference2}) instead of (\ref{eq:density-difference}).
We omit the details.
\end{pf}

\begin{remark}
We observe that $K$ and $\alpha$ defined in the above theorem
are independent of $t$. Therefore, under
Assumption \ref{assum}(i)--(iii), we have $\alpha(k) \leq K\alpha^{k}$
($\alpha$-mixing, defined in (\ref{eq:strong-mix}))
and under Assumption \ref{assum}(i)--(ii) and (iv), $\beta(k) \leq
K\alpha^{k}$
($\beta$-mixing, defined in (\ref{eq:betamix})) for all
$\sqrt{1-\delta}<\alpha<1$.

Moreover, since $\sigma(X_{t+k})\subset\sigma(X_{t+k},\ldots,X_{t+p-1})$
and $\sigma(X_{t})\subset\sigma(X_{t},\ldots,X_{t-p+1})$,
the 2-mixing rate is also geometric with $\tilde{\alpha}(k)\leq
K\alpha^{k}$
($\tilde{\alpha}(k)$ defined in (\ref{eq:2-mix})).
\end{remark}

\section{Mixing for $\operatorname{ARCH}(\infty)$ processes}\label{sec:ARCHinfty}

In this section, we derive mixing rates for the $\operatorname{ARCH}(\infty)$
process. We first define the process and state the assumptions that we
will use.

\subsection{The $\operatorname{ARCH}(\infty)$ process}

The $\operatorname{ARCH}(\infty)$ process has many interesting features,
which are useful in several applications.
For example, under certain conditions on the coefficients,
the $\operatorname{ARCH}(\infty)$ process can exhibit `near long memory' behaviour
(see \cite{pgir-00}).
The squares of the $\operatorname{ARCH}(\infty)$ process satisfy the representation
\begin{eqnarray}
\label{eq:ARCHinfty}
X_{t} = Z_{t}\Biggl(a_{0} + \sum_{j=1}^{\infty}a_{j}X_{t-j} \Biggr),
\end{eqnarray}
where $Z_{t}$ are i.i.d.~positive random variables with $\mathbb
{E}(Z_{t}) = 1$ and
$a_{j}$ are positive parameters. With a slight abuse of terminology, we
will call
the squared $\operatorname{ARCH}(\infty)$ process an $\operatorname{ARCH}(\infty)$ process. It is
worth mentioning that
the $\operatorname{GARCH}(p,q)$ process has an $\operatorname{ARCH}(\infty)$ representation, where the
$a_{j}$ decay geometrically with $j$. Giraitis and Robinson \cite{pgir-rob-01}, Robinson and Zaffaroni \cite
{pob-06} and
Subba~Rao \cite{psub-06} consider parameter estimation for the $\operatorname{ARCH}(\infty)$
process.

We will use Assumption \ref{assum} and the assumptions below.

\begin{assumption}\label{assum2}
\textup{(i)} We have $\sum_{j=1}^{\infty}a_{j} <
1-\delta$ and $a_{0}> 0 $.\vspace*{-6pt}
\begin{itemize}
\item[(ii)] For some $\nu>1$,
$\mathbb{E}|X_{t}|^{\nu}<\infty$
(we note that this is fulfilled if
$[\mathbb{E}|Z_{0}^{\nu}|]^{1/\nu}\sum_{j=1}^{\infty}a_{j}<1$).
\end{itemize}
\end{assumption}

Giraitis \textit{et al.} \cite{pgir-00} have shown that
under Assumption \ref{assum2}(i),
the $\operatorname{ARCH}(\infty)$ process has a stationary solution and a finite mean
(i.e., $\sup_{t\in\mathbb{Z}}\mathbb{E}(X_{t}) < \infty$).
It is worth mentioning that
since the $\operatorname{ARCH}(\infty)$ process has a stationary solution,
the shift $t$ plays no role when obtaining mixing bounds, that is,
$\sup_{G\in\sigma(X_{k+t}),H\in\sigma(X_{t})}|P(G\cap H)-P(G)P(H)|=
\sup_{G\in\sigma(X_{k}),H\in\sigma(X_{0})}|P(G\cap H)-P(G)P(H)|$.
Furthermore, the conditional density of
$X_{t+k}$ given $\underline{Z}_{t+k-1}^{t+1}$ and $\underline
{X}^{-\infty}_{\,t}$
is not a function of $t$. Hence, in the section below, we let
$f_{0,k}$ denote the conditional density of
$X_{t+k}$ given $(\underline{Z}_{t+k-1}^{t+1}$ and
$\underline{X}^{-\infty}_{\,t})$ and
for $s\geq1$, let
$f_{s,k}$ denote the conditional density of $X_{t+k+s}$
given $(\underline{X}_{\,t+k+s-1}^{t+k},
\underline{Z}_{t+k-1}^{t}$ and $\underline{X}^{-\infty}_{\,t})$.


\subsection{The $\operatorname{ARCH}(\infty)$ process and the Volterra series expansion}

We now write $X_{k}$ in terms of $\underline{Z}_{k-1}^{1}$ and
$\underline{X} = (X_{0},X_{-1},\ldots)$ and use this to derive
the conditional densities $f_{0,k}$ and $f_{s,k}$.
It can be seen from the result below (equation (\ref{eq:Xk}))
that, in general, the
$\operatorname{ARCH}(\infty)$ process is not Markovian.

\begin{lemma}
\label{lemma:ARCHinftyexpan}
Suppose that $\{X_{t}\}$ satisfies (\ref{eq:ARCHinfty}).
Then
\begin{eqnarray}
\label{eq:Xk}
X_{k} = \mathcal{P}_{0,k}(\underline{Z})Z_{k} +
\mathcal{Q}_{0,k}(\underline{Z},\underline{X})Z_{k},
\end{eqnarray}\vspace*{-18pt}

\noindent
where
\begin{eqnarray}\label{eq:PQ}
\mathcal{P}_{0,k}(\underline{Z}) &=& \Biggl[a_{0}+
\sum_{m=1}^{k}\sum_{k=j_{m} > \cdots>j_{1}>0}
\Biggl(\prod_{i=1}^{m-1}a_{j_{i+1} - j_{i}}\Biggr)
\Biggl( \prod_{i=1}^{m-1}Z_{j_{i}}\Biggr)
\Biggr], \nonumber
\\[-8pt]\\[-8pt]
\mathcal{Q}_{0,k}(\underline{Z},\underline{X}) &=&
\sum_{r=1}^{k}\Biggl\{
\sum_{m=1}^{k}\sum_{k=j_{m} > \cdots>j_{1}=r}
\Biggl( \prod_{i=1}^{m-1}a_{j_{i+1} - j_{i}}\Biggr)
\Biggl( \prod_{i=1}^{m-1}Z_{j_{i}}\Biggr)\Biggr\}
\,{d}_{r}(\underline{X}).\nonumber
\end{eqnarray}
Furthermore, setting $\mathcal{Q}_{0,k} = 0$
for $k \geq1$, we have that
$\mathcal{Q}_{0,k}(\underline{Z},\underline{X})$
satisfies the recursion
$\mathcal{Q}_{0,k}(\underline{Z},\underline{X}) =
\sum_{j=1}^{k}a_{j}\mathcal{Q}_{0,k-j}
(\underline{Z},\underline{X})Z_{k-j}
+ d_{k}(\underline{X})$,
where $d_{k}(\underline{X}) =
\sum_{j=0}^{\infty}a_{k+j}X_{-j}$ (for $k\geq1$).
\end{lemma}

\begin{pf} This can be found in Appendix A.3 of the technical
report.
\end{pf}

We will use the result above to derive the 2-mixing rate.
To derive $\alpha$ and $\beta$ mixing, we require the
density of $X_{k+s}$ given $\underline{X}_{\, k+s-1}^{k}$,
$\underline{Z}_{k-1}^{1}$
and $\underline{X}_{\,0}^{-\infty}$, which uses the following lemma.

\begin{lemma}
\label{lemma:4.2}
Suppose that $\{X_{t}\}$ satisfies (\ref{eq:ARCHinfty}).
For $s\geq1$, we then have
\begin{eqnarray}\label{eq:Xks-expan}
X_{k+s} &=& Z_{k+s}\{\mathcal{P}_{s,k}(\underline{Z})+
\mathcal{Q}_{s,k}(\underline{Z},\underline{X})\nonumber
\},
\\
&&\hspace*{-20pt}\mbox{where }
\mathcal{P}_{s,k}(\underline{Z})
= a_{0} +
\sum_{j=1}^{s}a_{j}X_{k+s-j} +
\sum_{j=s+1}^{\infty}a_{j}Z_{k+s-j}\mathcal{P}_{0,k+s-j}
(\underline{Z}),
\\
\mathcal{Q}_{s,k}(\underline{Z},\underline{X}) &=&
\sum_{j=s+1}^{k+s}a_{j}Z_{k+s-j}\mathcal{Q}_{0,k+s-j}(\underline{Z},
\underline{X}) + d_{k+s}(\underline{X}). \nonumber
\end{eqnarray}
\end{lemma}

\begin{pf} This can be found in Appendix A.3 of the technical
report.
\end{pf}

Using (\ref{eq:Xk}) and (\ref{eq:Xks-expan}), for all $s\geq0$, we
have that
$Z_{k+s} = \frac{X_{k+s}}{
\mathcal{P}_{s,k}(\underline{Z})+
\mathcal{Q}_{s,k}(\underline{Z},\underline{X})}$,
which leads to the conditional densities
\begin{eqnarray}\label{eq:PQbar}
&& f_{s,k}(y_{s}|\underline{y}{}_{s-1}, \underline{z},
\underline{x}) =\frac{1}{\mathcal{P}_{s,k}(\underline{z})+
\mathcal{Q}_{s,k}(\underline{z},\underline{x})}
f_{Z}\biggl(\frac{y_{s}}{\mathcal{P}_{s,k}(\underline{z})+
\mathcal{Q}_{s,k}(\underline{z},\underline{x})} \biggr).
\end{eqnarray}
In the proofs below, $\mathcal{Q}_{0,k}(\underline
{1}_{\,k-1},\underline{x})$
plays a prominent role.
By using the recursion in
Lemma \ref{lemma:ARCHinftyexpan} and (\ref{eq:PQbar}), setting
$\underline{x} = \underline{X}_{\,0}^{-\infty}$ and noting that
$\mathbb{E}(\mathcal{Q}_{s,k}(\underline{Z},\underline{x}))=
\mathcal{Q}_{s,k}(\underline{1}_{\,k-1},\underline{x})$,
we obtain the recursion
$\mathcal{Q}_{0,k}(\underline{1}_{\,k-1},\underline{x}) =
\sum_{j=1}^{k}a_{j+s}\mathcal{Q}_{0,k-j}(\underline
{1}_{\,k-j-1},\underline{x})
+ d_{k+s}(\underline{x})$. We use this to obtain a solution for
$\mathcal{Q}_{0,k}(\underline{1}_{\,k-1},\underline{x})$ in terms of
$\{d_{k}(\underline{x})\}_{k}$ in the lemma below.

\begin{lemma}\label{lemma:MAinfty}
Suppose that $\{X_{t}\}$ satisfies (\ref{eq:ARCHinfty}) and
Assumption \ref{assum2} is fulfilled.
There then exists $\{\psi_{j}\}$ such that for all $|z|\leq1$, we have
$(1-\sum_{j=1}^{\infty}a_{j}z^{j})^{-1} = \sum_{j=0}^{\infty}\psi
_{j}z^{j}$.
Furthermore, if $\sum_{j}|j^{\alpha}a_{j}| < \infty$, then
\cite{phan-86} have shown that $\sum_{j}|j^{\alpha}\psi_{j}| <
\infty$.
For $k\leq0$, set $d_{k}(\underline{x}) = 0$ and
$\mathcal{Q}_{0,k}(\underline{1}_{\,k-1},\underline{x}) = 0$. For
$k\geq1$,
$\mathcal{Q}_{0,k}(\underline{1}_{\,k-1},\underline{x})$ then has the solution\vspace*{-3pt}
\begin{eqnarray}
\label{eq:MAinfty1}
\mathcal{Q}_{0,k}(\underline{1}_{\,k-1},\underline{x}) =
\sum_{j=0}^{\infty}\psi_{j}d_{k-j}(\underline{x}) =
\sum_{j=0}^{k-1}\psi_{j}d_{k-j}(\underline{x})
= \sum_{j=0}^{k-1}
\psi_{j}\Biggl\{\sum_{i=0}^{\infty}a_{k-j+i}x_{-i}\Biggr\},\vspace*{-3pt}
\end{eqnarray}
where $\underline{x} = (x_{0},x_{-1},\ldots)$.
\end{lemma}

\begin{pf} This appears in Appendix A.3 of the technical report. \vspace*{-2pt}
\end{pf}

\subsection{Mixing for $\operatorname{ARCH}(\infty)$ processes}\vspace*{-2pt}

In this section, we show that the mixing rates are not necessarily
geometric and depend on the rate of
decay of the coefficients $\{a_{j}\}$ (we illustrate this in the
following example).
Furthermore, for $\operatorname{ARCH}(\infty)$ processes, the strong
mixing rate and 2-mixing rate can be different.

\begin{example}
Let us consider the $\operatorname{ARCH}(\infty)$ process, $\{X_{t}\}$,
defined in (\ref{eq:ARCHinfty}).
Giraitis \textit{et al.} \cite{pgir-00} have shown that if
$a_{j}\sim j^{-(1+\delta)}$ (for some $\delta>0$) and
$[\mathbb{E}(Z_{t}^{2})]^{1/2}\sum_{j=1}^{\infty}a_{j} <1$,
then $|\operatorname{cov}(X_{0},X_{k})|\sim k^{-(1+\delta)}$.
That is, the absolute sum of the covariances is finite, but `only just' if
$\delta$ is small.
If $Z_{t} < 1$, it is straightforward to see that
$X_{t}$ is a bounded random variable and by using Ibragimov's
inequality (see \cite{bhal-hey-80}), we have\vspace*{-3pt}
\begin{eqnarray*}
|\operatorname{cov}(X_{0},X_{k})| \leq C
\sup_{A\in\sigma(X_{0}),B\in\sigma(X_{k})} |P(A\cap B) - P(A)P(B)|\vspace*{-3pt}
\end{eqnarray*}
for some $C<\infty$. Noting that
$|\operatorname{cov}(X_{0},X_{k})|= \mathrm{O}(k^{-(1+\delta)})$, this gives
a lower bound of $\mathrm{O}(k^{-(1+\delta)})$ on the 2-mixing rate.
\end{example}

To obtain the mixing rates we will use Proposition
\ref{propgeneral-bound}, this result requires bounds on
$\mathcal{D}_{s,k} = |f_{s,k}(y_{s}|\underline{y}{}_{s-1},
\underline{z},\underline{x})-f_{s,k}(y_{s}|\underline{y}{}_{s-1},
\underline{z},0)|$ and its integral.

\begin{lemma}\label{lemma:H2}
Suppose that $\{X_{t}\}$ satisfies (\ref{eq:ARCHinfty})
and let $\mathcal{D}_{s,k}$ and $\mathcal{Q}_{0,k}(\cdot)$
be defined as in (\ref{eq:G}) and (\ref{eq:PQ}), respectively.
If Assumptions \textup{\ref{assum}(iii)} and \ref{assum2} are fulfilled, then\vspace*{-2pt}
\begin{eqnarray}
\label{eq:lemma:eq1}
&&\mathbb{E}_{\underline{Z}}\biggl(
\int|f_{0,k}(y|\underline{Z},\underline{x}) -
f_{0,k}(y|\underline{Z},0)|\,\mathrm{d}y \biggr)\nonumber
\\[-9pt]\\[-9pt]
 &&\quad\leq
\frac{\mathcal{Q}_{0,k}(\underline{1}_{\,k-1},\underline{x})}{a_{0}}
=\sum_{j=0}^{k-1}
|\psi_{j}|\Biggl\{\sum_{i=0}^{\infty}a_{k-j+i}x_{-i}\Biggr\}\nonumber
\end{eqnarray}

\noindent
and, for $s\geq1$,
\begin{eqnarray}
\label{eq:lemma:eq3}
&&\mathbb{E}_{\underline{Z}}\biggl(\sup_{\underline{y}{}_{s-1}\in
\mathbb{R}^{s}}\int\mathcal{D}_{s,k}
(y_{s}|\underline{y}{}_{s-1},\underline{Z},
\underline{x})\,\mathrm{d}y_{s}\biggr)\nonumber
\\[-8pt]\\[-8pt]
&&\quad \leq
\frac{1}{a_{0}}\Biggl\{
\sum_{j=s+1}^{k+s}a_{j} \sum_{l=0}^{k+s-j}|\psi_{l}|
\sum_{i=0}^{\infty}a_{k+s-j-l+i}x_{-i}
+ \sum_{i=0}^{\infty}a_{k+s+i}x_{-i}\Biggr\}. \nonumber
\end{eqnarray}
If Assumptions \textup{\ref{assum}(iv)} and \ref{assum2} are fulfilled
and $\mathcal{E}$ is defined as in (\ref{eq:defE}), then
\begin{eqnarray}
\label{eq:lemma:eq3beta}
&&\mathbb{E}_{\underline{Z}}\biggl(\sup_{\underline{y}{}_{s-1}\in
\mathbb{R}^{s}}\int\sup_{\underline{x}\in\mathcal{E}}\mathcal{D}_{s,k}
(y_{s}|\underline{y}{}_{s-1},\underline{Z},
\underline{x})\,\mathrm{d}y_{s} \biggr) \nonumber
\\[-8pt]\\[-8pt]
&&\quad \leq
\frac{1}{a_{0}}\Biggl\{
\sum_{j=s+1}^{k+s}a_{j} \sum_{l=0}^{k+s-j}|\psi_{l}|
\sum_{i=0}^{\infty}a_{k+s-j-l+i}\eta_{-i}
+ \sum_{i=0}^{\infty}a_{k+s+i}\eta_{-i}\Biggr\},\nonumber
\end{eqnarray}
where $\underline{x}=(x_{0},x_{-1},\ldots)$ is a positive vector.
\end{lemma}

\begin{pf} This can be found in Appendix A.3 of the technical
report.
\end{pf}


We require the following simple lemma to prove the theorem below.

\begin{lemma}\label{lemma:bias-variance-min}
If $\{c_{i}\}$, $\{d_{i}\}$ and $\{\eta_{-i}\}$
are positive sequences, then\vspace*{-2pt}
\begin{eqnarray}\label{eq:min}
\inf_{\underline{\eta}}
\Biggl\{ \sum_{i=0}^{\infty}(c_{i}\eta_{-i} +
d_{i}\eta_{-i}^{-\nu} )\Biggr\}
= \bigl(\nu^{{1/(1+\nu)}} + \nu^{-{\nu/(\nu+1)}}\bigr)
\sum_{i=0}^{\infty}c_{i}^{{\nu/(\nu+1)}}d_{i}^{{1/(\nu+1)}}.\vspace*{-2pt}
\end{eqnarray}
\end{lemma}

\begin{pf} This appears in Appendix A.3 of the technical report.
\end{pf}

In the following theorem, we obtain $\alpha$-mixing and $\beta$-mixing
bounds for the $\operatorname{ARCH}(\infty)$ process.

\begin{theorem}\label{thm:infty}
Suppose that $\{X_{t}\}$ satisfies (\ref{eq:ARCHinfty}).
\begin{enumerate}[(a)]
\item[(a)] Suppose Assumptions \textup{\ref{assum}(iii)}
and \ref{assum2} hold. We then have
\begin{eqnarray}\label{eq:lower-boundinfty}
&&\sup_{G\in\mathcal{F}^{k}_{\infty}, H\in
\mathcal{F}_{0}^{-\infty}}
|P(G\cap H) - P(G)P(H)| \nonumber
\\[-2pt]
&&\quad \leq K(\nu)\sum_{i=0}^{\infty}
\Biggl[\frac{1}{a_{0}}\sum_{s=0}^{\infty}
\sum_{j=s+1}^{k+s}a_{j}\sum_{l=0}^{k+s-j}|\psi_{l}|
a_{k+s-j-l+i}
\\[-2pt]
&&{}\qquad\hspace*{41pt} + \frac{1}{a_{0}}\sum_{s=0}^{\infty}
a_{k+s+i}\Biggr]^{{\nu/(\nu+1)}}
[\mathbb{E}|X_{0}|^{\nu}]^{{1/(\nu+1)}}
,\nonumber
\end{eqnarray}
where $K(\nu) = 3(\nu^{{1/(1+\nu)}} +
\nu^{-{\nu/(\nu+1)}})$.
\begin{enumerate}[(ii)]
\item[(i)] If the parameters of the $\operatorname{ARCH}(\infty)$ process satisfy
$|a_{j}| \sim j^{-\delta}$ and $|\psi_{j}|\sim j^{-\delta}$ ($\psi
_{j}$ defined in
Lemma \ref{lemma:MAinfty}), then we have
\begin{eqnarray*}
\sup_{G\in\mathcal{F}^{k}_{\infty},H\in
\mathcal{F}_{0}^{-\infty}}|P(G\cap H) - P(G)P(H)| \leq
K\cdot[ k (k+1)^{-\tilde{\delta}+3}+
(k+1)^{-\tilde{\delta}+2}],
\end{eqnarray*}
where $\tilde{\delta} = \delta\times(\frac{\nu}{\nu+1})$.

\item[(ii)] If the parameters of the $\operatorname{ARCH}(\infty)$ process satisfy
$|a_{j}| \sim\delta^{j}$ and $\psi_{j}\sim\delta^{j}$, where
$0< \delta< 1$
(an example is the $\operatorname{GARCH}(p,q)$ process), then we have
\begin{eqnarray*}
\sup_{G\in\mathcal{F}^{k}_{\infty},H\in\mathcal
{F}^{0}_{-\infty}}
|P(G\cap H) - P(G)P(H)|
&\leq& C \cdot k\cdot\delta^{k/2},
\end{eqnarray*}
where $C$ is a finite constant.
\end{enumerate}
\item[(b)]
If Assumptions \textup{\ref{assum}(iv)} and \ref{assum2} hold, then we have
\begin{eqnarray}\label{eq:lower-boundinfty3}
&&\sup_{\{G_{i}\}\in\mathcal{F}^{k}_{\infty},
\{H_{j}\}\in\mathcal{F}_{0}^{-\infty}}
\sum_{i}\sum_{j}
|P(G_{i}\cap H_{j}) - P(G_{i})P(H_{j})| \nonumber
\\
&&\quad \leq K(\nu)\sum_{i=0}^{\infty}
\Biggl[\frac{1}{a_{0}}\sum_{s=0}^{\infty}
\sum_{j=s+1}^{k+s}a_{j}\sum_{l=0}^{k+s-j}|\psi_{l}|
a_{k+s-j-l+i}
\\
&&{}\qquad\hspace*{41pt} + \frac{1}{a_{0}}\sum_{s=0}^{\infty}
a_{k+s+i}\Biggr]^{{\nu/(\nu+1)}}
[\mathbb{E}|X_{0}|^{\nu}]^{{1/(\nu+1)}},\nonumber
\end{eqnarray}
where $\{G_{i}\}$ and $\{H_{j}\}$ are partitions of $\Omega$.
We mention that the bounds for the $\alpha$-mixing rates
for different orders of $\{a_{j}\}$ and $\{\psi_{j}\}$ derived in
(\textup{i}) also hold for the $\beta$-mixing rate.
\end{enumerate}
\end{theorem}

\begin{pf} We first prove (a). We use the fact that
\begin{eqnarray*}
\sup_{G\in\mathcal{F}^{k}_{\infty},H
\in\mathcal{F}_{0}^{-\infty}}
|P(G\cap H) - P(G)P(H)|=
\lim_{n\rightarrow\infty}
\sup_{G\in
\mathcal{F}^{k}_{k+n}}{H\in\mathcal{F}_{0}^{-\infty}}
|P(G\cap H) - P(G)P(H)|
\end{eqnarray*}
and find a bound for each $n$.
By using (\ref{eq:mix-bd}) to bound
$\sup_{G\in
\mathcal{F}^{k}_{k+n},H\in\mathcal{F}_{0}^{-\infty}}
|P(G\cap H) - P(G)P(H)|$,
we see that for all sets $\mathcal{E}$
(as defined in (\ref{eq:defE})), we have
\begin{eqnarray}\label{eq:GH1}
&&\sup_{G
\in\mathcal{F}^{k}_{k+n},H\in
\underline{F}_{\,0}^{-\infty}}
|P(G\cap H) - P(G)P(H)|\nonumber
\\
&&\quad \leq
2\sup_{\underline{x}\in\mathcal{E}}
\sum_{s=0}^{n}
\mathbb{E}_{\underline{Z}}
\biggl(\sup_{\underline{y}{}_{s-1}\in\mathbb{R}^{s}}
\biggl\{\int\mathcal{D}_{s,k}(y_{s}|
\underline{y}{}_{s-1},\underline{Z},\underline{x})\,\mathrm{d}y_{s}\biggr\}
\biggr)
\\
&&{}\qquad +
4P(X_{0}>\eta_{0}\mbox{ or, } \ldots,
X_{-n}> \eta_{-n}).\nonumber
\end{eqnarray}
To bound the integral in (\ref{eq:GH1}), we use
(\ref{eq:lemma:eq3}) to obtain
\begin{eqnarray*}
&&\sup_{\underline{x}\in\mathcal{E}}
\sum_{s=0}^{n}
\mathbb{E}_{\underline{Z}}\biggl(
\sup_{\underline{y}{}_{s-1}\in
\mathbb{R}^{s}}\int_{\mathbb{R}}
\mathcal{D}_{s,k}(y_{s}|
\underline{y}{}_{s-1},\underline{Z},\underline{x})\,\mathrm{d}y_{s} \biggr)
\nonumber
\\
&&\quad =\frac{1}{a_{0}}\sum_{s=0}^{n}\Biggl\{
\sum_{j=s+1}^{k+s}a_{j} \sum_{l=0}^{k+s-j}|\psi_{l}|
\sum_{i=0}^{\infty}a_{k+s-j-l+i}\eta_{-i} +
\sum_{i=0}^{\infty}a_{k+s+i}\eta_{-i}\Biggr\}.
\end{eqnarray*}
Now, by using Markov's inequality, we have that
$P(X_{0}>\eta_{0}\mbox{ or},\ldots,X_{-n}\geq\eta_{-n}) \leq\break
\sum_{i=0}^{n}\frac{\mathbb{E}(|X_{i}|^{\nu})}{\eta_{-i}^{\nu}}$.
Substituting this and the above
into (\ref{eq:GH1}) and letting $n\rightarrow\infty$ gives
\begin{eqnarray}
\label{eq:strong-mixing}\vspace{-3mm}
&&\sup_{G\in\mathcal{F}^{k}_{\infty},H\in
\mathcal{F}^{-\infty}_{0}}
|P(G\cap H) - P(G)P(H)| \nonumber
\\
&&\quad \leq\inf_{\eta}\Biggl[
\frac{2}{a_{0}}\sum_{s=0}^{\infty}\Biggl\{
\sum_{j=s+1}^{k+s}a_{j} \sum_{l=0}^{k+s-j}|\psi_{l}|
\sum_{i=0}^{\infty}a_{k+s-j-l+i}\eta_{-i}
\\
&&{}\qquad\hspace*{29pt}\hspace*{19pt} +
\sum_{i=0}^{\infty}a_{k+s+i}\eta_{-i}\Biggr\}
+ 4\mathbb{E}|X_{0}|^{\nu}
\sum_{i=0}^{\infty}\eta_{-i}^{-\nu}\Biggr],\nonumber
\end{eqnarray}
where $\eta= (\eta_{0},\eta_{-1},\ldots)$.

We now use (\ref{eq:min}) to minimize (\ref{eq:strong-mixing}), which
gives us
(\ref{eq:lower-boundinfty}). The proof of (i) can be found in the
technical report.
It is straightforward to prove (ii) using (\ref{eq:min}).

The proof of (b) is very similar to the proof of (a), but uses
(\ref{eq:lemma:eq3beta})
rather than (\ref{eq:lemma:eq3}). We omit the details.
\end{pf}

\begin{remark}
Under the assumptions of Theorem \textup{\ref{thm:infty}(a)}, we have a bound
for the $\alpha$-mixing rate, that is,
$\alpha(k)\leq\zeta(k)$, where
$\zeta(k)=
K[\frac{1}{a_{0}}\sum_{s=0}^{\infty}
\sum_{j=s+1}^{k+s}a_{j}\sum_{l=0}^{k+s-j}|\psi_{l}|
a_{k+s-j-l+i} + \frac{1}{a_{0}}\sum_{s=0}^{\infty}
a_{k+s+i}]^{{\nu/(\nu+1)}}$.
Under the assumptions of Theorem \textup{\ref{thm:infty}(a)}, the
$\beta$-mixing coefficient is bounded by
$\beta(k)\leq\zeta(k)$.
\end{remark}

In the following theorem, we consider a bound for the 2-mixing rate of an
$\operatorname{ARCH}(\infty)$ process.

\begin{theorem}\label{thm:infty2}
Suppose that $\{X_{t}\}$ satisfies (\ref{eq:ARCHinfty}) and that
Assumption \textup{\ref{assum}(iii)} and \ref{assum2}
hold. We then have
\begin{eqnarray}\label{eq:lower-boundinfty2}
&&\sup_{G\in\sigma(X_{k}),H\in\mathcal{F}_{0}^{-\infty}}
|P(G\cap H) - P(G)P(H)|\nonumber
\\[-8pt]\\[-8pt]
&&\quad \leq K(\nu)\sum_{i=0}^{\infty}
\Biggl[\frac{1}{a_{0}}
\sum_{j=0}^{k-1}a_{j}|\psi_{j}|
a_{k-j+i}\Biggr]^{{\nu/(\nu+1)}}
[\mathbb{E}|X_{0}|^{\nu}]^{{1/(\nu+1)}},\nonumber
\end{eqnarray}
where $K(\nu) = 3(\nu^{{1/(1+\nu)}} + \nu^{-{\nu/(\nu+1)}})$.

If the parameters of the $\operatorname{ARCH}(\infty)$ process satisfy
$a_{j} \sim j^{-\delta}$ and $|\psi_{j}|\sim j^{-\delta}$ ($\psi
_{j}$ defined in
Lemma \ref{lemma:MAinfty}), then we have
\begin{eqnarray}\label{eq:2mix-poly}
\sup_{G\in\sigma(X_{k}),H\in
\mathcal{F}_{0}^{-\infty}}|P(G\cap H) - P(G)P(H)|
\leq K \cdot k(k+1)^{-\tilde{\delta}+1},
\end{eqnarray}
where $\tilde{\delta} = \delta\times(\frac{\nu}{\nu+1})$.
\end{theorem}

\begin{pf} We use a similar proof to that of Theorem \ref
{thm:infty}. The
integral difference is replaced with the bound in (\ref{eq:lemma:eq1}) and
we again use Markov's inequality: together they give the bound
\begin{eqnarray}
\label{eq:Mmixing}
&&\sup_{G\in\sigma(X_{k}),H \in
\mathcal{F}_{0}^{-\infty}}|P(G\cap H) - P(G)P(H)|\nonumber
\\[-8pt]\\[-8pt]
&&\quad \leq \inf_{\underline{\eta}}\Biggl[2
\frac{1}{a_{0}}\sum_{j=0}^{k-1}
|\psi_{j}|\Biggl\{\sum_{i=0}^{\infty}a_{k-j+i}\eta_{-i}\Biggr\}
+ 4\mathbb{E}|X_{0}|^{\nu}
\sum_{i=0}^{\infty}\frac{1}{\eta_{-i}^{\nu}}\Biggr].\nonumber
\end{eqnarray}
Finally, to obtain (\ref{eq:lower-boundinfty2}) and (\ref{eq:2mix-poly}),
we use (\ref{eq:Mmixing}) and a proof similar to that of
Theorem \ref{thm:infty}(i).  We omit the details.
\end{pf}

\begin{remark}
Comparing (\ref{eq:2mix-poly}) and Theorem \ref{thm:infty}(i), we see that
the 2-mixing bound is of a smaller order than the strong mixing bound.

In fact, it could well be that the 2-mixing bound is of a smaller order than
Theorem \ref{thm:infty2}(i). This is because Theorem \ref{thm:infty2}(i)
gives a bound for $\sup_{G\in\sigma(X_{k}),H\in
\sigma(X_{0},X_{-1},\ldots)}
|P(G\cap H) - P(G)P(H)|$, whereas the 2-mixing bound
restricts the $\sigma$-algebra of the left tail to $\sigma(X_{0})$.
However, we have not been able to show this and this is a problem
that requires further consideration.
\end{remark}



\begin{appendix}
\section*{Appendix: Proofs}\label{app}

\renewcommand{\thesubsection}{A.\arabic{subsection}}
\setcounter{subsection}{0}
\subsection{\texorpdfstring{Proof of
Proposition \protect\ref{propgeneral-bound}}{Proof of Proposition 2.1}}\label{sec:app1}

\renewcommand{\theequation}{\arabic{equation}}
\setcounter{equation}{35}

\renewcommand{\thelemma}{A.\arabic{lemma}}
\setcounter{lemma}{0}

We will use the following three lemmas to prove
Proposition \ref{propgeneral-bound}.
\begin{lemma}\label{lemma:A1}
Let $G\in\mathcal{F}_{t+k+r_{2}}^{t+k}=
\sigma(\underline{X}_{\,t+k+r_{2}}^{t+k})$ and
$H, E \in\mathcal{F}_{t}^{t-r_{1}}= \sigma(\underline{X}_{\,t}^{t-r_{1}})$
(where $E$ is defined in (\ref{eq:defE})), and
use the notation of Proposition \ref{propgeneral-bound}. We then have
\begin{eqnarray}
\label{eq:A1a}
&& |P(G\cap H \cap E) - P(G\cap E)P(H)| \nonumber
\\
&&\quad \leq 2P(H) \sup_{\underline{x}\in
\mathcal{E}}\bigl|P(G|\underline{X}^{t-r_{1}}_{\,t} = \underline{x}) -
P(G|\underline{X}^{t-r_{1}}_{\,t} = 0)\bigr|
\\
&&\qquad{} +\inf_{\underline{x}\in\mathcal{E}}P(G|\underline
{X}^{t-r_{1}}_{t} =
\underline{x})\{P(H)P(E^{c}) + P(H\cap E^{c}) \}.\nonumber
\end{eqnarray}
%
\end{lemma}

\begin{pf} To prove the result, we first observe that
\begin{eqnarray*}
P(G\cap H \cap E) &=& P\bigl(\underline{X}_{\,t+k+r_{2}}^{t+k}\in
\mathcal{G},
\underline{X}_{\,t}^{t-r_{1}}\in(\mathcal{H}\cap\mathcal{E})\bigr)
\\
&=&
\int_{\mathcal{H}\cap\mathcal{E}}\int_{\mathcal{G}}
\mathrm{d}P(\underline{X}_{\,t}^{t-r_{1}}\leq\underline{y},\underline
{X}_{t+k+r_{2}}^{t+k}
\leq\underline{x}) \nonumber
\\
&=& \int_{\mathcal{H}\cap\mathcal{E}} \biggl\{\int_{\mathcal{G}}
\mathrm{d}P(\underline{X}_{\,t+k+r_{2}}^{t+k}\leq\underline{y}|
\underline{X}_{\,t}^{t-r_{1}} = \underline{x})
\biggr\}\,\mathrm{d}P(\underline{X}_{\,t}^{t-r_{1}}\leq\underline{x}) \nonumber
\\
&=&
\int_{\mathcal{H}\cap\mathcal{E}} P(\underline
{X}_{t+k+r_{2}}^{t+k}\in
\mathcal{G}|\underline{X}_{\,t}^{t-r_{1}} = \underline{x})
\,\mathrm{d}P(\underline{X}_{\,t}^{t-r_{1}}\leq\underline{x}).
\end{eqnarray*}
Therefore, by using the above and the fact that $P(H\cap E)\leq P(H)$,
we obtain the following inequalities:
\begin{eqnarray}\label{eq:50}
&&\inf_{\underline{x}\in\mathcal{E}}
P(\underline{X}_{\,t+k+r_{2}}^{t+k}\in
\mathcal{G}|\underline{X}_{\,t}^{t-r_{1}} = \underline{x})
P(H\cap E)\nonumber
\\[-9pt]\\[-9pt]
&&\quad \leq
P(G\cap H \cap E) \leq\sup_{\underline{x}\in\mathcal{E}}
P(\underline{X}_{\,t+k+r_{2}}^{t+k}\in
\mathcal{G}|\underline{X}_{\,t}^{t-r_{1}} = \underline
{x})P(H)\nonumber
\end{eqnarray}
and
\begin{eqnarray}\label{eq:51}
&&\inf_{\underline{x}\in\mathcal{E}}
P(\underline{X}_{\,t+k+r_{2}}^{t+k}\in
\mathcal{G}|\underline{X}_{\,t}^{t-r_{1}} = \underline{x})
P(E)\nonumber
\\[-9pt]\\[-9pt]
&&\quad \leq
P(G \cap E) \leq\sup_{\underline{x}\in\mathcal{E}}
P(\underline{X}_{\,t+k+r_{2}}^{t+k}\in
\mathcal{G}|\underline{X}_{\,t}^{t-r_{1}} = \underline{x})P(E).
\nonumber
\end{eqnarray}
Subtracting (\ref{eq:50}) from (\ref{eq:51}) and using
$P(H\cap E) = P(H) - P(H\cap E^{c})$ gives the inequalities
\begin{eqnarray}\label{eq:GHEmore}
&&P(G\cap H \cap E) - P(G\cap E)P(H)\nonumber
\\
&&\quad \leq
\sup_{\underline{x}\in\mathcal{E}}
P(\underline{X}_{\,t+k+r_{2}}^{t+k}\in
\mathcal{G}|\underline{X}_{\,t}^{t-r_{1}} = \underline{x})P(H)
\\[-2pt]
&&\qquad {} - \inf_{\underline{x}\in\mathcal{E}}
P(\underline{X}_{\,t+k+r_{2}}^{t+k}\in
\mathcal{G}|\underline{X}_{\,t}^{t-r_{1}} = \underline{x})
P(H) + P(E^{c})P(H),\nonumber
\\\label{eq:GHEless}
&&P(G\cap H \cap E) - P(G\cap E)P(H)\nonumber
\\[-2pt]
&&\quad \geq
\inf_{\underline{x}\in\mathcal{E}}
P(\underline{X}_{\,t+k+r_{2}}^{t+k}\in
\mathcal{G}|\underline{X}_{\,t}^{t-r_{1}} = \underline{x})P(H)
\\[-2pt]
&&{}\qquad  - \sup_{\underline{x}\in\mathcal{E}}
P(\underline{X}_{\,t+k+r_{2}}^{t+k}\in
\mathcal{G}|\underline{X}_{\,t}^{t-r_{1}} = \underline{x})
P(H) - P(E^{c}\cap H).\nonumber
\end{eqnarray}
Combining (\ref{eq:GHEmore}) and (\ref{eq:GHEless}), we obtain
\begin{eqnarray}
\label{eq:A1b}
&&|P(G\cap H \cap E) - P(G\cap E)P(H)| \nonumber
\\[-2pt]
&&\quad \leq P(H)\Bigl| \sup_{\underline{x}\in
\mathcal{E}}P(G|\underline{X}^{t-r_{1}}_{\,t} = \underline{x}) -
\inf_{\underline{x}\in
\mathcal{E}}P(G|\underline{X}^{t-r_{1}}_{\,t} = \underline{x})\Bigr|
\\[-2pt]
&&{}\qquad  +\inf_{\underline{x}\in\mathcal{E}}P(G|\underline
{X}^{t-r_{1}}_{t} =
\underline{x})\{P(H)P(E^{c}) + P(H\cap E^{c}) \}.\nonumber
\end{eqnarray}
Using the triangle inequality, we have
\begin{eqnarray*}
\Bigl| \sup_{\underline{x}\in
\mathcal{E}}P(G|\underline{X}^{t-r_{1}}_{\,t} = \underline{x}) -
\inf_{\underline{x}\in
\mathcal{E}}P(G|\underline{X}^{t-r_{1}}_{\,t} = \underline{x})\Bigr|
\leq 2\sup_{\underline{x}\in
\mathcal{E}}| P(G|\underline{X}^{t-r_{1}}_{\,t} = \underline{x}) -
P(G|\underline{X}^{t-r_{1}}_{\,t} = 0)|.
\end{eqnarray*}
Substituting the above into (\ref{eq:A1b}) gives
(\ref{eq:A1a}), as required.
\end{pf}

We now obtain a bound for the first term on the right-hand
side of (\ref{eq:A1a}).

\begin{lemma}\label{lemma:A2}
Let $f_{\underline{X}_{\,t+k+r_{2}}^{t+k}|\underline{X}_{\,t}^{t-r_{1}}}$ denote
the density of $\underline{X}_{\,t+k+r_{2}}^{t+k}$ given
$\underline{X}_{\,t}^{t-r_{1}}$ and
$\mathcal{G}$ and $\mathcal{H}$ be defined as in (\ref{eq:HG}). Then,
\begin{eqnarray}\label{eq:A21}
\bigl| P(G|\underline{X}^{t-r_{1}}_{\,t} = \underline{x}) -
P(G|\underline{X}^{t-r_{1}}_{\,t} = 0)\bigr|
\leq\int_{\mathcal{G}}\mathcal{D}_{0,k,t}
(\underline{y}|\underline{x})\, \mathrm{d}\underline{y}.
\end{eqnarray}
Let $\underline{W}_{t+k-1}^{t+1}$ be a random vector which is independent
of $\underline{X}_{\,t}^{t-r_{1}}$ and let $f_{\underline{W}}$
denote the density of $\underline{W}_{t+k-1}^{t+1}$. If
$G\in\sigma(X_{t+k})$, then
\begin{eqnarray}\label{eq:A22}
\int_{\mathcal{G}}
\bigl|f_{X_{t+k}|\underline{X}_{\,t}^{t-r_{1}}}
(y|\underline{x}) -
f_{X_{t+k}|\underline{X}_{\,t}^{t-r_{1}}}(y|0)
\bigr|\,\mathrm{d}y \leq\mathbb{E}_{\underline{W}}
\biggl(\int_{\mathbb{R}}\mathcal{D}_{0,k,t}(y|\underline
{W},\underline{x})
\,\mathrm{d}y \biggr)
\end{eqnarray}
and if $G\in\sigma(\underline{X}_{\,t+k+r_{2}}^{t+k})$, then
\begin{eqnarray} \label{eq:A23}
&&\int_{\mathcal{G}}
\bigl|f_{\underline{X}_{\,t+k+r_{2}}^{t+k}|\underline{X}_{\,t}^{t-r_{1}}}
(\underline{y}|\underline{x}) -
f_{\underline{X}_{\,t+k+r_{2}}^{t+k}|\underline
{X}_{t}^{t-r_{1}}}(\underline{y}|0)
\bigr|\,\mathrm{d}\underline{y}\nonumber
\\[-8pt]\\[-8pt]
&&\quad \leq \sum_{s=0}^{r_{2}}\mathbb{E}_{\underline{W}}
\biggl(\sup_{\underline{y}{}_{s-1}}\int_{{\bf G}_{s}}
\mathcal{D}_{s,k,t}(y_{s}|\underline{y}{}_{s-1},\underline
{w},\underline{x})
\,\mathrm{d}y_{s} \biggr).\nonumber
\end{eqnarray}
%
\end{lemma}

\begin{pf} The proof of (\ref{eq:A21}) is clear from the
definition of $\mathcal{D}_{s,k,t}$,
hence we omit the details.

To prove (\ref{eq:A22}), we first note that by independence
of $\underline{W}^{t+1}_{t+k-1}$ and $\underline{X}_{\,t}^{t-r_{2}}$,
we have
that $f_{\underline{W}|\underline{X}_{\,t}^{t-r_{1}}}(
\underline{w}|\underline{x})=f_{\underline{W}}(\underline{w})
$, where $f_{\underline{W}|\underline{X}_{\,t}^{t-r_{1}}}$
is the conditional density of $\underline{W}^{t+1}_{t+k-1}$ given
$\underline{X}_{\,t}^{t-r_{1}}$. Therefore, we have
\begin{eqnarray*}
f_{X_{t+k}|\underline{X}_{\,t}^{t-r_{1}}}(y|\underline{x})
&=&
\int_{\mathbb{R}^{k-1}} f_{X_{t+k}|\underline{W},
\underline{X}_{\,t}^{t-r_{1}}} (y|\underline{w},\underline{x})
f_{\underline{W}}(\underline{w})\,\mathrm{d}\underline{w} =
\int_{\mathbb{R}^{k-1}} f_{0,k,t} (y|\underline{w},\underline{x})
f_{\underline{W}}(\underline{w})\,\mathrm{d}\underline{w}.
\end{eqnarray*}
Substituting the above into
$\int_{\mathcal{G}}|f_{X_{t+k}|\underline{X}_{\,t}^{t-r_{1}}}
(y|\underline{x}) -f_{X_{t+k}|\underline{X}_{\,t}^{t-r_{1}}}(y|0)|\,\mathrm{d}y$
and using the definition of $\mathbb{E}_{\underline{W}}$
now gives (\ref{eq:A22}).

To prove (\ref{eq:A23}), we note that, by using the same argument used
to prove
(\ref{eq:A22}), we have
\begin{eqnarray}
\label{eq:cond-den}
f_{\underline{X}_{\,t+k+r_{2}}^{t+k}|
\underline{X}_{\,t}^{t-r_{1}}}(\underline{y}|\underline{x}) =
\int_{\mathbb{R}^{k-1}} f_{\underline{W}}(\underline{w})\prod_{s=0}^{r_{2}}
f_{s,k,t}(y_{s}|\underline{y}{}_{s-1},\underline{w},
\underline{x})\, \mathrm{d}\underline{w}.
\end{eqnarray}
Now, repeatedly subtracting and adding $f_{s,k,t}$ gives
\begin{eqnarray}\label{eq:A2expan}
&& f_{\underline{X}_{\,t+k+r_{2}}^{t+k}|
\underline{X}_{\,t}^{t-r_{1}}}(\underline{y}|\underline{x}) -
f_{\underline{X}_{\,t+k+r_{2}}^{t+k}|
\underline{X}_{\,t}^{t-r_{1}}}(\underline{y}| 0 )\nonumber
\\
&&\quad =\sum_{s=0}^{r_{2}}\int_{\mathbb{R}^{k-1}}
f_{\underline{W}}(\underline{w})\Biggl\{\prod_{a=0}^{s-1}
f_{a,k,t}(y_{a}|\underline{y}{}_{a-1},\underline{w},
\underline{x})\Biggr\}\nonumber
\\[-8pt]\\[-8pt]
&&\qquad{}\hspace*{37pt}\times\Biggl\{\prod
_{b=s+1}^{r_{2}}f_{b,k,t}(y_{b}|\underline{y}{}_{b-1},
\underline{w},0)\Biggr\}\nonumber
\\
&&{}\qquad\hspace*{37pt}\times \{f_{s,k,t}(y_{s}|\underline{y}{}_{s-1},\underline{w},
\underline{x}) - f_{s,k,t}(y_{s}|\underline{y}{}_{s-1},\underline
{w},0) \}
\,\mathrm{d}\underline{w}. \nonumber
\end{eqnarray}
Therefore, taking the integral of the above over $\mathcal{G}$ gives
\begin{eqnarray}\label{eq:A2expan2}
&& \int_{\mathcal{G}} \bigl|f_{\underline{X}_{\,t+k+r_{2}}^{t+k}|
\underline{X}_{\,t}^{t-r_{1}}}(\underline{y}|\underline{x}) -
f_{\underline{X}_{\,t+k+r_{2}}^{t+k}|
\underline{X}_{\,t}^{t-r_{1}}}(\underline{y}| 0 )
\bigr|\,\mathrm{d}\underline{y}\nonumber
\\
&&\quad \leq \sum_{s=0}^{r_{2}}
\int_{\mathbb{R}^{k-1}}
f_{\underline{W}}(\underline{w})\Biggl\{\Biggl[
\prod_{a=0}^{s-1}\int_{\mathcal{G}_{a}}
f_{a,k,t}(y_{a}|\underline{y}{}_{a-1},\underline{w},\underline
{x})\,\mathrm{d}y_{a}\nonumber
\\[-8pt]\\[-8pt]
&&{}\qquad\hspace*{72pt}\times \prod_{b=s+1}^{r_{2}}\int_{\mathcal{G}_{b}}
f_{b,k,t}(y_{b}|\underline{y}{}_{b-1},\underline{w},\underline{x})
\,\mathrm{d}y_{b}\Biggr]\nonumber
\\
&&{}\qquad\hspace*{72pt} \times\sup_{\underline{y}{}_{s-1}}
\int_{\mathcal{G}_{s}}
\bigl|f_{s,k,t}(y_{s}|\underline{y}{}_{s-1},\underline{w},
\underline{x}) - f_{s,k,t}(y_{s}|\underline{y}{}_{s-1},\underline
{w},0) \bigr|
\,\mathrm{d}y_{s}\Biggr\} \,\mathrm{d}\underline{w}.\quad \nonumber
\end{eqnarray}
Next, we observe that since $\mathcal{G}_{j}\subset
\mathbb{R}$ and $\int_{\mathbb{R}}f_{s,k,t}(y_{s}|\underline{y}{}_{s-1},
\underline{w},\underline{x})\,\mathrm{d}y_{s}=1$,
we\vspace*{-3pt} have\break
$(\prod_{a=0}^{s-1}\int_{\mathcal{G}_{a}}
f_{a,k,t}(y_{a}|\underline{y}{}_{a-1},\underline{w},\underline{x})
\,\mathrm{d}y_{a})
(\prod_{b=s+1}^{r_{2}}\int_{\mathcal{G}_{b}}
f_{b,k,t}(y_{b}|\underline{y}{}_{b-1},\underline{w},\underline{x})
\,\mathrm{d}y_{b})\leq1$.
Finally, substituting this bound into
(\ref{eq:A2expan2}) gives (\ref{eq:A23}).
\end{pf}

The following lemma will be used to show $\beta$-mixing and
uses the above lemmas.

\begin{lemma}
Suppose that $\{G_{i}\}\in\mathcal{F}_{t+k+r_{2}}^{t+k}$,
$\{H_{j}\}\in\mathcal{F}_{t}^{t-r_{1}}$ and
$\{G_{i}\}$ and $\{H_{j}\}$ are partitions of~$\Omega$. We then have
\begin{eqnarray}
\label{eq:A31}
&&\sum_{i,j}|P(G_{i}\cap H_{j}\cap E ) - P(G_{i}\cap E)P(H_{j})|
\nonumber
\\[-8pt]\\[-8pt]
&&\quad \leq 2\sum_{i}\sup_{\underline{x}\in\mathcal{E}}
\bigl|P(G_{i}|\underline{X}_{\,t}^{t-r_{1}} = \underline{x}) -
P(G_{i}|\underline{X}_{\,t}^{t-r_{1}} = 0)\bigr| + 2P(E^{c})\quad \mbox{and
}\nonumber\\[-20pt]\nonumber
\end{eqnarray}
\begin{eqnarray}
\label{eq:A32}
\sum_{i,j}|P(G_{i}\cap H_{j}\cap E^{c}) - P(G_{i}\cap E^{c})P(H_{j})|
\leq2P(E^{c}).
\end{eqnarray}
\end{lemma}

\begin{pf} Substituting the inequality in (\ref{eq:A1a}) into
$\sum_{i,j}|P(G_{i}\cap H_{j}\cap E ) - P(G_{i}\cap E)P(H_{j})|$ gives
\begin{eqnarray}\label{eq:60a}
&&\sum_{i,j}|P(G_{i}\cap H_{j}\cap E ) - P(G_{i}\cap E)P(H_{j})|
\nonumber
\\[-2pt]
&&\quad \leq 2\sum_{j} P(H_{j})\sum_{i} \sup_{\underline{x}\in
\mathcal{E}}\bigl|P(G_{i}|\underline{X}^{\,t-r_{1}}_{t} = \underline
{x}) -
P(G_{i}|\underline{X}^{t-r_{1}}_{\,t} = 0)\bigr|
\\[-2pt]
&&\qquad{}+ \sum_{i,j}
\inf_{\underline{x}\in\mathcal{E}}P(G_{i}|\underline
{X}^{t-r_{1}}_{t} =
\underline{x})\{P(H_{j})P(E^{c}) + P(H_{j}\cap E^{c}) \}.\nonumber
\end{eqnarray}
The sets $\{H_{j}\}$ are partitions of $\Omega$, hence
$\sum_{i}P(H_{j}) = 1$ and $\sum_{i}P(H_{j}\cap E^{c})\leq1$.
Using these observations together with (\ref{eq:60a}) gives (\ref{eq:A31}).

Inequality (\ref{eq:A32}) immediately follows from the fact that
$\{H_{j}\}$ and $\{G_{i}\}$ are disjoint sets.
\end{pf}

Using the above three lemmas, we can now prove Proposition
\ref{propgeneral-bound}.

\begin{pf*}{Proof of Proposition \ref{propgeneral-bound}, equation (\ref
{eq:mix-bd})}
It is straightforward to show that
\begin{eqnarray*}
| P(G\cap H) - P(G)P(H)| &\leq& | P(G\cap H \cap E) -
P(G\cap
E)P(H)|
\\[-2pt]
&&{}+ | P(G\cap H \cap E^{c}) - P(G\cap E^{c})P(H)|.
\end{eqnarray*}
Now, by substituting (\ref{eq:A21}) into (\ref{eq:A1a}) and using the
above, we get
\begin{eqnarray*}
 | P(G\cap H) - P(G)P(H)| &\leq&
2\sup_{\underline{x}\in\mathcal{E}}\int_{\mathcal{G}}
\bigl|f_{\underline{X}_{\,t+k+r_{2}}^{t+k}|\underline{X}_{\,t}^{t-r_{1}}}
(\underline{y}|\underline{x}) -
f_{\underline{X}_{\,t+k+r_{2}}^{t+k}|\underline
{X}_{t}^{t-r_{1}}}(\underline{y}|0)
\bigr|\,\mathrm{d}\underline{y}
\\[-2pt]
&&{} +\inf_{\underline{x}\in\mathcal{E}}P(G|\underline
{X}^{t-r_{1}}_{t} =
\underline{x})\{P(H)P(E^{c}) + P(H\cap E^{c}) \}
\\[-2pt]
&&{}+
P(G\cap H \cap E^{c}) + P(G \cap E^{c})P(H).
\end{eqnarray*}
Finally, by using the facts that $\mathcal{G}\subset\mathbb{R}^{r_{2}+1}$,
$P(G\cap H \cap E^{c})\leq P(E^{c})$, $P(G \cap E^{c})P(H)\leq P(E^{c})$
and $\inf_{\underline{x}\in\mathcal{E}}P(G|\underline
{X}^{t-r_{1}}_{t} =
\underline{x})\leq1$, we obtain (\ref{eq:mix-bd}).
\end{pf*}

\begin{pf*}{Proof of Proposition \ref{propgeneral-bound}, equation
(\ref{eq:mix-bdii})}
It is worth noting that the
proof of (\ref{eq:mix-bdii}) is similar to the proof of (\ref{eq:mix-bd}).
Using (\ref{eq:A31}) and
the same arguments as those in the proof of (\ref{eq:mix-bd}), we have
\begin{eqnarray}\label{eq:mix-bdiiP}
&&\sum_{i,j}|P(G_{i}\cap H_{j}) - P(G_{i})P(H_{j})|\nonumber
\\[-2pt]
&&\quad \leq 2\sum_{i}\sup_{\underline{x}\in\mathcal{E}}
\int_{\mathcal{G}_{i}}
\bigl|f_{\underline{X}_{\,t+k+r_{2}}^{t+k}|\underline{X}_{\,t}^{t-r_{1}}}
(\underline{y}|\underline{x}) -
f_{\underline{X}_{\,t+k+r_{2}}^{t+k}|\underline
{X}_{t}^{t-r_{1}}}(\underline{y}|0)
\bigr|\,\mathrm{d}\underline{y} + 4P(E^{c}) \nonumber
\\[-9pt]\\[-9pt]
&&\quad \leq 2\sum_{i}
\int_{\mathcal{G}_{i}}\sup_{\underline{x}\in\mathcal{E}}
\bigl|f_{\underline{X}_{\,t+k+r_{2}}^{t+k}|\underline{X}_{\,t}^{t-r_{1}}}
(\underline{y}|\underline{x}) -
f_{\underline{X}_{\,t+k+r_{2}}^{t+k}|\underline
{X}_{t}^{t-r_{1}}}(\underline{y}|0)
\bigr|\,\mathrm{d}\underline{y} + 4P(E^{c}) \nonumber
\\[-2pt]
&&\quad \leq 2
\int_{\mathbb{R}^{r_{2}+1}}\sup_{\underline{x}\in\mathcal{E}}
\bigl|f_{\underline{X}_{\,t+k+r_{2}}^{t+k}|\underline{X}_{\,t}^{t-r_{1}}}
(\underline{y}|\underline{x}) -
f_{\underline{X}_{\,t+k+r_{2}}^{t+k}|\underline
{X}_{t}^{t-r_{1}}}(\underline{y}|0)
\bigr|\,\mathrm{d}\underline{y} + 4P(E^{c}),\nonumber
\end{eqnarray}
where $H_{j} =
\{\omega; \underline{X}_{\,t}^{t-r_{1}}(\omega) \in\mathcal{H}_{j}\}
$ and
$G_{i} = \{\omega; \underline{X}_{\,t+k+r_{2}}^{t+k}(\omega)
\in\mathcal{G}_{i}\}$, which gives (\ref{eq:mix-bdii}).
\end{pf*}

\begin{pf*}{Proof of Proposition \ref{propgeneral-bound}, equation
(\ref{eq:mix-bd2})}
To prove the result, we substitute the bound in (\ref{eq:A23})
into (\ref{eq:mix-bd}) to obtain (\ref{eq:mix-bd2}).
\end{pf*}

\begin{pf*}{Proof of Proposition \ref{propgeneral-bound}, equation
(\ref{eq:mix-bd2ii})} To prove (\ref{eq:mix-bd2ii}), we substitute
(\ref{eq:A23}) into (\ref{eq:mix-bdii}) to obtain (\ref{eq:mix-bd2ii}).
\end{pf*}

\subsection{\texorpdfstring{Proofs in Section \protect\ref{sec:tvARCH}}{Proofs in Section 3}}\label{sec:app2}

\begin{pf*}{Proof of Lemma \ref{lemma:3.1}}
We first prove (\ref{eq:tvARCHEXPa}) with $s=0$. Suppose that $k\geq
1$. Focusing on the first element of $\underline{X}_{\,t+k}^{t+k-p+1}$ in
(\ref{eq:underlineX}) and factoring out $Z_{t+k}$ gives
\begin{eqnarray*}
X_{t+k} &=& Z_{t+k}\Biggl\{a_{0}(t+k) +
 \Biggl[\tilde{A}_{t+k} \sum_{r=0}^{k-2}
\prod_{i=1}^{r}A_{t+k-i}(Z)b_{t+k-r-1}(Z)\Biggr]_{1}
\\
&&{}\hspace*{24pt} +
\Biggl[\tilde{A}_{t+k}
\Biggl\{\prod_{i=1}^{k-1}A_{t+k-i}(Z)\Biggr\}\underline{X}_{\,t}^{t-p+1}\Biggr]_{1}
\Biggr\},
\end{eqnarray*}
which is (\ref{eq:tvARCHEXPa}) (with $s=0$).
To prove (\ref{eq:tvARCHEXPa}) for $1\leq s\leq p$, we note that using the
tvARCH$(p)$ representation in (\ref{eq:tvARCH})
and (\ref{eq:tvARCHEXPa}) for $s=0$ gives
\begin{eqnarray*}
X_{t+k+s} &=&
Z_{t+k+s}\Biggl\{a_{0}(t+k+s) + \sum_{i=1}^{s-1}a_{i}(t+k+s)X_{t+k+s-i}
+ \sum_{i=s}^{p}a_{i}(t+k+s)X_{t+k+s-i}
\Biggr\} \nonumber
\\
&=& Z_{t+k+s}\{
\mathcal{P}_{s,k,t}(
\underline{Z})
+ \mathcal{Q}_{s,k,t}(\underline{Z},\underline{X})\},
\end{eqnarray*}
where $\mathcal{P}_{s,k,t}$ and $\mathcal{Q}_{s,k,t}$ are defined in
(\ref{eq:Q}). Hence, this gives (\ref{eq:tvARCHEXPa}).
Since $a_{j}(\cdot)$ and $Z_{t}$ are positive, it is clear that
$\mathcal{P}_{s,k,t}$ and $\mathcal{Q}_{s,k,t}$ are positive
random variables.
\end{pf*}

\begin{pf*}{Proof of Lemma \ref{lemma:smaller-algebra}}
We first note that since $\{X_{t}\}$ satisfies a tvARCH$(p)$
representation ($p<\infty$) it is $p$-Markovian, hence for any $r_{2}>p$,
the $\sigma$-algebras generated\vspace*{-2pt} by $\underline{X}_{\,t+k+r_{2}}^{t+k}$ and
$(\underline{Z}_{t+k+r_{2}}^{t+k+p},\underline{X}_{\,t+k+p-1}^{t+k})$
are the same. Moreover, by using the fact that for all $\tau> t$,
$Z_{\tau}$ is
independent of $X_{t}$, we have
%
\begin{eqnarray}
&&\sup_{G\in\mathcal{F}^{t+k}_{\infty},H\in
\mathcal{F}^{-\infty}_{t}}
|P(G\cap H) - P(G)P(H)|\nonumber
\\[-8pt]\\[-8pt]
&&\quad =
\sup_{G\in\mathcal{F}_{t+k+p-1}^{t+k},H\in
\mathcal{F}^{t-p+1}_{t}}|P(G\cap H) - P(G)P(H)|.\nonumber
\end{eqnarray}
Now, by using the above,
Proposition \ref{propgeneral-bound}, equation (\ref{eq:mix-bd2}),
and the fact that $\underline{Z}_{t+k-1}^{t+1}$ and
$\underline{X}_{\,t}^{t-p+1}$ are independent, for any set
$\mathcal{E}$ (defined as in (\ref{eq:defE})), we have
\begin{eqnarray}
&&
\sup_{G\in\mathcal{F}_{t+k+p-1}^{t+k},H\in
\mathcal{F}^{t-p+1}_{t}}|P(G\cap H) - P(G)P(H)| \nonumber
\\
&&\quad \leq 2 \sup_{\underline{x}\in\mathcal{E}}
\sum_{s=0}^{p-1} \mathbb{E}_{\underline{Z}}\biggl(\sup_{\underline
{y}_{s-1}
\in\mathbb{R}^{s}}\int
\mathcal{D}_{s,k,t}(y_{s}|\underline{y}{}_{s-1},\underline
{z},\underline{x})
\,\mathrm{d}y_{s}\biggr)
\\
&&{}\qquad + 4P(X_{t}>\eta_{0}\mbox{ or}, \ldots,
X_{t-p+1}>\eta_{-p+1}).\qquad\nonumber
\end{eqnarray}
Finally, using the fact that $P(X_{t}>\eta_{0}\mbox{ or
}X_{t-1}>\eta_{-1}, \ldots,
X_{t-p+1}>\eta_{-p+1}) \leq\break \sum_{j=0}^{p-1} P(X_{t-j} > \eta_{-j})$
gives (\ref{eq:diff-f}).

The proof of (\ref{eq:diff-f2}) is similar to the proof above,
but uses (\ref{eq:mix-bd2ii}) instead of (\ref{eq:mix-bd2}), so
we omit the details.
\end{pf*}

We require the following simple lemma to prove
Lemmas \ref{lemma:H} and \ref{lemma:H2}.

\begin{lemma}
If Assumption \textup{\ref{assum}(iii)} is satisfied, then, for any positive
$A$ and $B$, we have
\begin{eqnarray}
\label{eq:AB1}
\int_{\mathbb{R}}
\biggl|\frac{1}{A+B}f_{Z}\biggl(\frac{y}{A+B}\biggr) - \frac{1}{A}f_{Z}\biggl(\frac{y}{A}\biggr)\biggr|\,\mathrm{d}y
\leq K\biggl(\frac{B}{A} + \frac{B}{A+B}\biggr).
\end{eqnarray}
If Assumption \textup{\ref{assum}(iv)} is satisfied, then, for any
positive $A$, positive continuous function $B\dvtx \mathbb
{R}^{r_{2}+1}\rightarrow
\mathbb{R}$ and set $E$ (defined as in (\ref{eq:defE})), we have
\begin{eqnarray}
\label{eq:AB12}
\int_{\mathbb{R}}\sup_{\underline{x}\in\mathcal{E}}
\biggl|\frac{1}{A+B(\underline{x})}f_{Z}\biggl(\frac{y}{A+B(\underline{x})}\biggr) -
\frac{1}{A}f_{Z}\biggl(\frac{y}{A}\biggr)\biggr|\,\mathrm{d}y
\leq K\sup_{\underline{x}\in\mathcal{E}}\biggl(\frac{B(\underline
{x})}{A} +
\frac{B(\underline{x})}{A+B(\underline{x})}\biggr).\quad
\end{eqnarray}
\end{lemma}

\begin{pf} To prove (\ref{eq:AB1}), we observe that
\begin{eqnarray*}
\int_{\mathbb{R}}
\biggl|\frac{1}{A+B}f_{Z}\biggl(\frac{y}{A+B}\biggr) - \frac{1}{A}f_{Z}\biggl(\frac
{y}{A}\biggr)\biggr|\,\mathrm{d}y = I + \mathit{II},
\end{eqnarray*}
where
\begin{eqnarray*}
 I&=&
\int_{\mathbb{R}} \frac{1}{A+B}\biggl|f_{Z}\biggl(\frac{y}{A+B}\biggr) - f_{Z}\biggl(\frac
{y}{A}\biggr)\biggr|\,\mathrm{d}y
\quad\mbox{and}\quad
\mathit{II} = \int_{\mathbb{R}}\biggl(\frac{1}{A+B} - \frac{1}{A}\biggr)f_{Z}\biggl(\frac{y}{A}\biggr).
\end{eqnarray*}
To bound $I$, we note that by changing variables with $u = y/(A+B)$
and under Assumption~\ref{assum}(iii), we get
\begin{eqnarray*}
I \leq\int_{\mathbb{R}}
\biggl|f_{Z}(u) - f_{Z}\biggl(u\biggl(1+\frac{B}{A}\biggr)\biggr)
\biggr|\,\mathrm{d}u\leq
K\frac{B}{A}.
\end{eqnarray*}
It is straightforward to show that $\mathit{II}\leq\frac{B}{A+B}$. Hence, the bounds
for $I$ and $\mathit{II}$ give (\ref{eq:AB1}).

The proof of (\ref{eq:AB12}) is the same as above,
but uses Assumption~\ref{assum}$(iv)$ instead of Assumption~\ref
{assum}(iii),
so we omit the details.
\end{pf}

\begin{pf*}{Proof of Lemma \ref{lemma:H}} We first show that
\begin{eqnarray}
\label{eq:kappa-bound}
\sup_{\underline{y}{}_{s-1}\in\mathbb{R}^{s}}\int
\mathcal{D}_{s,k,t}(y_{s}|
\underline{y}{}_{s-1},\underline{z},\underline{x})\,\mathrm{d}y_{s}
\leq \frac{K}{\inf_{t\in\mathbb{Z}}a_{0}(t)}
\mathcal{Q}_{s,k,t}(\underline{z},\underline{x})
\end{eqnarray}
and use this to prove (\ref{eq:density-difference}).
We note that when
$\underline{x} = 0$, $\mathcal{Q}_{s,k,t}(\underline{z},0) = 0$ and
$f_{s,k,t}(y_{s}|\underline{y}{}_{s-1},\underline{z},0)
= \mathcal{P}_{s,k,t}(\underline{z})^{-1}
f_{Z}(\frac{y_{s}}{\mathcal{P}_{s,k,t}
(\underline{z})})$. Therefore, using (\ref{eq:fxz}) gives
\begin{eqnarray*}
\mathcal{D}_{s,k,t}(y_{s}|
\underline{y}{}_{s-1},\underline{z},\underline{x})
&=& \biggl|\frac{1}{
\mathcal{P}_{s,k,t}(\underline{z}) +
\mathcal{Q}_{s,k,t}(\underline{z},\underline{x})}
f_{Z}\biggl(\frac{y_{s}}{\mathcal{P}_{s,k,t}(\underline{z})
+\mathcal{Q}_{s,k,t}(\underline{z},\underline{x})}\biggr)
\\
&&{}\ - \frac{1}{\mathcal{P}_{s,k,t}(\underline{x})}
f_{Z}\biggl(\frac{y_{s}}{\mathcal{P}_{s,k,t}(
\underline{z})}\biggr)\biggr|.
\end{eqnarray*}
Now, recalling that $\mathcal{P}_{s,k,t}$ and $\mathcal{Q}_{s,k,t}$
are both
positive and
setting $A = \mathcal{P}_{s,k,t}(\underline{z})$,
$B = \mathcal{Q}_{s,k,t}(\underline{z},\underline{x})$ and using
(\ref{eq:AB1}), we have
\begin{eqnarray*}
\int_{\mathbb{R}} \mathcal{D}_{s,k,t}(y_{s}|
\underline{y}{}_{s-1},\underline{z},\underline{x})\,\mathrm{d} y_{s} \leq
K\biggl( \frac{\mathcal{Q}_{s,k,t}(\underline{z},\underline{x})}{
\mathcal{P}_{s,k,t}(\underline{z})} +
\frac{\mathcal{Q}_{s,k,t}(\underline{z},\underline{x})}{
\mathcal{P}_{s,k,t}(\underline{z}) +
\mathcal{Q}_{s,k,t}(\underline{z},\underline{x})} \biggr).
\end{eqnarray*}
Finally, since $\mathcal{P}_{s,k,t}(
\underline{z}) > \inf_{t\in\mathbb{Z}}a_{0}(t)$, we have
$\int_{\mathbb{R}} \mathcal{D}_{s,k,t}(y_{s}|
\underline{y}{}_{s-1},\underline{z},\underline{x})\,\mathrm{d} y_{s} \leq
K\frac{\mathcal{Q}_{s,k,t}(\underline{z},\underline{x})}{
\inf_{t\in\mathbb{Z}}a_{0}(t)}$,
thus giving (\ref{eq:kappa-bound}).
By using (\ref{eq:kappa-bound}), we now prove (\ref{eq:density-difference}).
Substituting (\ref{eq:kappa-bound}) into the integral on the left-hand
side of (\ref{eq:density-difference}), using the fact that
$\mathbb{E}[\mathcal{Q}_{s,k,t}(\underline{Z},\underline{x})] =
\mathcal{Q}_{s,k,t}(\underline{1}_{\,k-1},\underline{x})$ and
substituting (\ref{eq:kappa-bound}) into (\ref{eq:diff-f}) gives
\begin{eqnarray}
\label{eq:68a}
\mathbb{E}_{\underline{Z}}\biggl( \sup_{\underline{y}{}_{s-1}\in
\mathbb{R}^{s}}
\int_{\mathbb{R}}\mathcal{D}_{s,k,t}(y_{s}|
\underline{y}{}_{s-1},\underline{Z},\underline{x})\,\mathrm{d}y_{s}\biggr)
\leq
K\frac{\mathbb{E}[\mathcal{Q}_{s,k,t}(\underline{Z},
\underline{x})]}{\inf_{t\in\mathbb{Z}}a_{0}(t)} =
K\frac{\mathcal{Q}_{s,k,t}(\underline{1}_{\,k-1},
\underline{x})}{\inf_{t\in\mathbb{Z}}a_{0}(t)}.\quad
\end{eqnarray}
We now find a bound for $\mathcal{Q}_{s,k,t}$.
By the definition of $\mathcal{Q}_{s,k,t}$ in (\ref{eq:Q}) and
using the matrix norm
inequality $[A\underline{x}]_{1}\leq K\|A\|_{\mathit{spec}}
\|\underline{x}\|$ ($\|\cdot\|_{\mathit{spec}}$ is the spectral norm), we have
\begin{eqnarray*}
\mathcal{Q}_{s,k,t}(\underline{1}_{\,k-1},\underline{x}) &=&
\sum_{i=s+1}^{p}a_{i}(t+k+s)\Biggl[A_{t+k+s-i}
\sum_{r=1}^{k+s-i}\Biggl\{\prod_{d=0}^{k+s-i}A_{t+k+s-i-d}\Biggr\}\underline
{x}\Biggr]_{1}
\\
&\leq& \frac{K}{\inf_{t\in\mathbb{Z}}a_{0}(t)}\sum_{i=s}^{p}a_{i}(t+k+s)
\Biggl\|A_{t+k+s-i}\Biggl\{\prod_{d=0}^{k-1}A_{t+k+s-i-d}\Biggr\}\Biggr\|_{\mathit{spec}}\|
\underline{x}\|.
\end{eqnarray*}
To bound the above, we note that
by Assumption \ref{assum}(i),
$\sup_{t\in\mathbb{Z}}
\sum_{j=1}^{p}a_{j}(t)\leq(1-\delta)$, therefore there exists a
$\tilde{\delta}$, where
$0<\tilde{\delta}< \delta < 1$ and such that, for all $t$, we have
$\|A_{t+k+s-i}\times\break \{\prod_{d=0}^{k-1}A_{t+k+s-i-d}\}\|_{\mathit{spec}}\leq
K(1-\tilde{\delta})^{k+1}$
for some finite $K$. Combining all of this gives
\begin{eqnarray}\label{eq:geoQ}
\mathcal{Q}_{s,k,t}(\underline{1}_{\,k-1},\underline{x}) &\leq&
\frac{K}{\inf_{t\in\mathbb{Z}}a_{0}(t)}\sum_{i=s}^{p}a_{i}(t+k+s)
\Biggl\|A_{t+k+s-i}\Biggl\{\prod_{d=0}^{k+s-i}A_{t+k+s-i-d}\Biggr\}\Biggr\|_{\mathit{spec}}
\|\underline{x}\| \nonumber
\\[-8pt]\\[-8pt]
&\leq& \frac{K}{\inf_{t\in\mathbb{Z}}a_{0}(t)}\sum_{i=s}^{p}a_{i}(t+k+s)
(1-\tilde{\delta})^{k+s-i}\|\underline{x}\|.\nonumber
\end{eqnarray}
Substituting the above into (\ref{eq:68a}) gives
(\ref{eq:density-difference}).

We now prove (\ref{eq:density-difference2}). We use the
same proof to show (\ref{eq:kappa-bound}), but apply (\ref{eq:AB1})
instead of (\ref{eq:AB12}) to obtain
%
\begin{eqnarray*}
\sup_{\underline{y}{}_{s-1}\in\mathbb{R}^{s}}\int
\sup_{\underline{x}\in\mathcal{E}}
\mathcal{D}_{s,k,t}(y_{s}|
\underline{y}{}_{s-1},\underline{z},\underline{x})\,\mathrm{d}y_{s}
\leq\frac{K}{\inf_{t\in\mathbb{Z}}a_{0}(t)}
\sup_{\underline{x}\in\mathcal{E}}
\mathcal{Q}_{s,k,t}(\underline{z},\underline{x}).\quad
\end{eqnarray*}
By substituting the above into (\ref{eq:diff-f2}) and using the
same proof to prove (\ref{eq:density-difference}), we obtain
\begin{eqnarray}\label{eq:68b}
&&\sum_{s=0}^{p-1} \int
\prod_{i=1}^{k-1}f_{Z}(z_{i})\sup_{\underline{y}{}_{s-1}\in
\mathbb{R}^{s}}\biggl\{\int_{\mathbb{R}} \sup_{\underline{x}\in
\mathcal{E}}
\mathcal{D}_{s,k,t}(y_{s}|
\underline{y}{}_{s-1},\underline{z},\underline{x})\,\mathrm{d}y_{s}\biggr\}
\,\mathrm{d}\underline{z}\nonumber
\\[-8pt]\\[-8pt]
&&\quad \leq K\frac{\mathbb{E}[\sup_{\underline{x}\in\mathcal{E}}
\mathcal{Q}_{s,k,t}(\underline{Z},
\underline{x})]}{\inf_{t\in\mathbb{Z}}a_{0}(t)}.\nonumber
\end{eqnarray}
Since $\mathcal{Q}_{s,k,t}(\underline{Z},
\underline{x})$ is a positive function and
$\sup_{\underline{x}\in\mathcal{E}}
\mathcal{Q}_{s,k,t}(\underline{Z},
\underline{x}) = \mathcal{Q}_{s,k,t}(\underline{Z},
\underline{\eta})$, we have
$\mathbb{E}[\sup_{\underline{x}\in\mathcal{E}}
\mathcal{Q}_{s,k,t}(\underline{Z},
\underline{x})]\leq\sup_{\underline{x}\in\mathcal{E}}
\mathbb{E}[\mathcal{Q}_{s,k,t}(\underline{Z},
\underline{x})] = \sup_{\underline{x}\in\mathcal{E}}
\mathcal{Q}_{s,k,t}(\underline{1}_{\,k-1},
\underline{x})$. Hence, by using (\ref{eq:geoQ}), we have
\begin{eqnarray*}
\frac{\mathbb{E}[\sup_{\underline{x}\in\mathcal{E}}
\mathcal{Q}_{s,k,t}(\underline{Z},
\underline{x})]}{\inf_{t\in\mathbb{Z}}a_{0}(t)} \leq
\frac{K(1-\tilde{\delta})^{k}\|\underline{x}\|}{\inf_{t\in\mathbb{Z}}
a_{0}(t)}.
\end{eqnarray*}
Substituting the above into (\ref{eq:68b})
gives (\ref{eq:density-difference2}).
\end{pf*}
\end{appendix}


\section*{Acknowledgements}

We would like to thank Piotr Kokoszka,
Mika Meitz and Joseph Tadjuidje
for several useful discussions. We also wish to thank the Associate
Editor and two anonymous referees for suggestions and comments which
greatly improved many aspects of the paper. The research of Suhasini
Subba Rao was partially
supported by an NSF Grant under DMS-0806096 and
the Deutsche Forschungsgemeinschaft under DA 187/15-1.

\printhistory

\end{document}